\documentclass[a4paper,10pt]{article}
\usepackage[right=1.35in, left=1.3in, footskip=.5in]{geometry}
\usepackage{amsmath,amsfonts,amsthm,wasysym}
\usepackage{mathrsfs}
\usepackage{graphicx}
\usepackage{float}
\usepackage{enumitem}
\usepackage{tikz}
\usetikzlibrary{calc}
\usepackage{multicol}
\usepackage{hyperref}

\theoremstyle{plain}
\newtheorem{Th}{Theorem}[section]
\newtheorem{Lem}[Th]{Lemma}
\newtheorem{Cor}[Th]{Corollary}
\newtheorem{Prop}[Th]{Proposition}

\theoremstyle{definition}
\newtheorem{Def}[Th]{Definition}
\newtheorem*{Not}{Notation}
\newtheorem{Ex}[Th]{Example}
\newtheorem{Assum}[Th]{Assumption}

\theoremstyle{remark}
\newtheorem*{Rem}{Remark}

\def\thm{\begin{Th}}
\def\endthm{\end{Th}}
\def\lemma{\begin{Lem}}
\def\endlemma{\end{Lem}}
\def\cor{\begin{Cor}}
\def\endcor{\end{Cor}}
\def\prop{\begin{Prop}}
\def\endprop{\end{Prop}}
\def\definition{\begin{Def}}
\def\enddefinition{\end{Def}}
\def\remark{\begin{Rem}}
\def\endremark{\end{Rem}}
\def\example{\begin{Ex}}
\def\endexample{\end{Ex}}
\def\demo{\begin{proof}}
\def\enddemo{\end{proof}}

\def\notation{\begin{Not}}
\def\endnotation{\end{Not}}
\def\assumption{\begin{Assum}}
\def\endassumption{\end{Assum}}

\def\BbR{\mathbb{R}}
\def\BbN{\mathbb{N}}

\def\e{\epsilon}
\def\a{\alpha}
\def\b{\beta}
\def\c{\gamma}

\def\vp{\varphi}

\def\s{\sigma}

\def\tF{\widetilde{\F}}

\def\word#1#2{{#1}_1\ldots{#1}_{#2}}

\def\H{\mathcal{H}}
\def\L{\mathcal{L}}
\def\E{\mathcal{E}}

\def\D{\mathcal{D}}

\def\S{\mathcal{S}}
\def\F{\mathcal{F}}
\def\G{\mathcal{G}}
\def\R{\mathcal{R}}

\def\SS{\Sigma}

\def\diam#1{{\rm diam}(#1)}

\def\ou{\overline{u}}

\def\RF{\mathcal{R}\mathcal{F}}

\def\hE{\widehat{\E}}
\def\hF{\widehat{\F}}

\def\sd#1#2{#1\backslash#2}

\def\l{\ell}

\def\MP{{\mathcal{M}\mathcal{P}}}

\title{Completely Symmetric Resistance Forms on the Stretched Sierpinski Gasket}
\author{P. Alonso Ruiz\thanks{Department of Mathematics, University of Connecticut, 341 Mansfield Rd, Unit-1009 Storrs, CT-06269; email: patricia.alonso-ruiz@uconn.edu}  \and U. Freiberg\thanks{Institute of Stochastics and Applications, University of Stuttgart, Pfaffenwaldring 57, 70569 Stuttgart, Germany; email: Uta.Freiberg@mathematik.uni-stuttgart.de} \and J. Kigami\thanks{Graduate School of Informatics, Kyoto University, Kyoto 606-8501, Japan; email: kigami@i.kyoto-u.ac.jp}}
\date{}

\begin{document}

\maketitle

\begin{abstract}
The stretched Sierpinski gasket, SSG for short, is the space obtained by replacing every branching point of the Sierpinski gasket by an interval. It has also been called the ``deformed Sierpinski gasket'' or ``Hanoi attractor''. As a result, it is the closure of a countable union of intervals and one might expect that a diffusion on SSG is essentially a kind of gluing of the Brownian motions on the intervals. In fact, there have been several works in this direction. There still remains, however, ``reminiscence'' of the Sierpinski gasket in the geometric structure of SSG and the same should therefore be expected for diffusions. This paper shows that this is the case. In this work, we identify all the completely symmetric resistance forms on SSG. A completely symmetric resistance form is a resistance form whose restriction to every contractive copy of SSG in itself is invariant under all geometrical symmetries of the copy, which constitute the symmetry group of the triangle.  We prove that completely symmetric resistance forms on SSG can be sums of the Dirichlet integrals on the intervals with some particular weights, or a linear combination of a resistance form of the former kind and the standard resistance form on the Sierpinski gasket.
\end{abstract}



\section{Introduction}
A major area of research interest in mathematical physics deals with the modelling of heat and wave propagation in branching media. One way to tackle this problem consists in approximating the object under consideration by unions of one-dimensional segments, and studying the combination of the corresponding equations on the segments. This approach has been extensively investigated under different names, for instance ``quantum graphs'' in mathematical physics~\cite{Kuc08} and ``cable systems'' in stochastic analysis~\cite{BB04}. 

\medskip

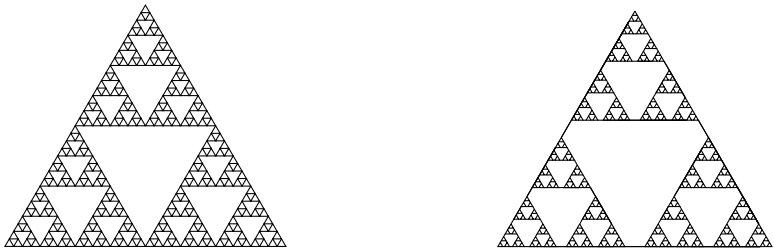
\begin{figure}[H]
\centering
\begin{tikzpicture}[scale=1/15]
\foreach \c in {0,1,2} {
\foreach \b in {0,1,2} {
\foreach \a in {0,1,2} {
\foreach \e in {0,1,2} {
\foreach \f in {0,1,2} {
\draw[] ($(90+120*\f:16)+(90+120*\e:8) + (90+120*\c:4)+(90+120*\b:2)+(90+120*\a:1)+(90:1)$) -- ($(90+120*\f:16)+(90+120*\e:8) + (90+120*\c:4)+(90+120*\b:2)+(90+120*\a:1)+(210:1)$) -- ($(90+120*\f:16)+(90+120*\e:8) + (90+120*\c:4)+(90+120*\b:2)+(90+120*\a:1)+(330: 1)$)--cycle; 
}
}
}
}
}
\end{tikzpicture}
\hspace*{1in}
\begin{tikzpicture}[scale=20/620]
\draw ($(90:2.75)$) --  ($(90+120:2.75)$) -- ($(90+120*2:2.75)$) -- ($(90:2.75)$) ($(90:2.75)+(120*2:2)$) to ++ ($(0:2)$) ($(90+120:2.75)+(0:2)$) to ++ ($(120:2)$)  ($(90+120*2:2.75)+(120:2)$) to ++ ($(120*2:2)$) ;
\draw ($(90:2.75)$) -- ($(90+120:2.75)+(60:6)$) ($(90+120*2:2.75)$) -- ($(90+120:2.75)+(0:6)$) ($(90:2.75)+(0:6)$) -- ($(90+120*2:2.75)+(60:6)$);
\foreach \a in {0,60} {
\draw ($(90:2.75)+(\a:6)$) --  ($(90+120:2.75)+(\a:6)$) -- ($(90+120*2:2.75)+(\a:6)$) -- ($(90:2.75)+(\a:6)$) ($(90:2.75)+(120*2:2)+(\a:6)$) to ++ ($(0:2)$) ($(90+120:2.75)+(0:2)+(\a:6)$) to ++ ($(120:2)$)  ($(90+120*2:2.75)+(120:2)+(\a:6)$) to ++ ($(120*2:2)$);
}
\draw ($(90:2.75)+(60:6)$) -- ($(90+120:2.75)+(60:6)+(60:13)$) ($(90+120*2:2.75)+(0:6)$) -- ($(90+120:2.75)+(0:6)+(0:13)$) ($(90:2.75)+(60:6)+(0:13)$) -- ($(90+120*2:2.75)+(0:6)+(60:13)$);
\foreach \b in {0,60} {
\draw ($(90:2.75)+(\b:13)$) --  ($(90+120:2.75)+(\b:13)$) -- ($(90+120*2:2.75)+(\b:13)$) -- ($(90:2.75)+(\b:13)$) ($(90:2.75)+(120*2:2)+(\b:13)$) to ++ ($(0:2)$) ($(90+120:2.75)+(0:2)+(\b:13)$) to ++ ($(120:2)$)  ($(90+120*2:2.75)+(120:2)+(\b:13)$) to ++ ($(120*2:2)$) ;
\draw ($(90:2.75)+(\b:13)$) -- ($(90+120:2.75)+(60:6)+(\b:13)$) ($(90+120*2:2.75)+(\b:13)$) -- ($(90+120:2.75)+(0:6)+(\b:13)$) ($(90:2.75)+(0:6)+(\b:13)$) -- ($(90+120*2:2.75)+(60:6)+(\b:13)$);
\foreach \a in {0,60} {
\draw ($(90:2.75)+(\a:6)+(\b:13)$) --  ($(90+120:2.75)+(\a:6)+(\b:13)$) -- ($(90+120*2:2.75)+(\a:6)+(\b:13)$) -- ($(90:2.75)+(\a:6)+(\b:13)$) ($(90:2.75)+(120*2:2)+(\a:6)+(\b:13)$) to ++ ($(0:2)$) ($(90+120:2.75)+(0:2)+(\a:6)+(\b:13)$) to ++ ($(120:2)$)  ($(90+120*2:2.75)+(120:2)+(\a:6)+(\b:13)$) to ++ ($(120*2:2)$) ;
}
}
\foreach \c in {0,60} {
\draw ($(90:2.75)+(\c:28)$) --  ($(90+120:2.75)+(\c:28)$) -- ($(90+120*2:2.75)+(\c:28)$) -- ($(90:2.75)+(\c:28)$) ($(90:2.75)+(120*2:2)+(\c:28)$) to ++ ($(0:2)$) ($(90+120:2.75)+(0:2)+(\c:28)$) to ++ ($(120:2)$)  ($(90+120*2:2.75)+(120:2)+(\c:28)$) to ++ ($(120*2:2)$) ;
\draw ($(90:2.75)+(\c:28)$) -- ($(90+120:2.75)+(60:6)+(\c:28)$) ($(90+120*2:2.75)+(\c:28)$) -- ($(90+120:2.75)+(0:6)+(\c:28)$) ($(90:2.75)+(0:6)+(\c:28)$) -- ($(90+120*2:2.75)+(60:6)+(\c:28)$);
\foreach \a in {0,60} {
\draw ($(90:2.75)+(\a:6)+(\c:28)$) --  ($(90+120:2.75)+(\a:6)+(\c:28)$) -- ($(90+120*2:2.75)+(\a:6)+(\c:28)$) -- ($(90:2.75)+(\a:6)+(\c:28)$) ($(90:2.75)+(120*2:2)+(\a:6)+(\c:28)$) to ++ ($(0:2)$) ($(90+120:2.75)+(0:2)+(\a:6)+(\c:28)$) to ++ ($(120:2)$)  ($(90+120*2:2.75)+(120:2)+(\a:6)+(\c:28)$) to ++ ($(120*2:2)$);
}
\draw ($(90:2.75)+(60:6)+(\c:28)$) -- ($(90+120:2.75)+(60:6)+(60:13)+(\c:28)$) ($(90+120*2:2.75)+(0:6)+(\c:28)$) -- ($(90+120:2.75)+(0:6)+(0:13)+(\c:28)$) ($(90:2.75)+(60:6)+(0:13)+(\c:28)$) -- ($(90+120*2:2.75)+(0:6)+(60:13)+(\c:28)$);
\foreach \b in {0,60} {
\draw ($(90:2.75)+(\b:13)+(\c:28)$) --  ($(90+120:2.75)+(\b:13)+(\c:28)$) -- ($(90+120*2:2.75)+(\b:13)+(\c:28)$) -- ($(90:2.75)+(\b:13)+(\c:28)$) ($(90:2.75)+(120*2:2)+(\b:13)+(\c:28)$) to ++ ($(0:2)$) ($(90+120:2.75)+(0:2)+(\b:13)+(\c:28)$) to ++ ($(120:2)$)  ($(90+120*2:2.75)+(120:2)+(\b:13)+(\c:28)$) to ++ ($(120*2:2)$) ;
\draw ($(90:2.75)+(\b:13)+(\c:28)$) -- ($(90+120:2.75)+(60:6)+(\b:13)+(\c:28)$) ($(90+120*2:2.75)+(\b:13)+(\c:28)$) -- ($(90+120:2.75)+(0:6)+(\b:13)+(\c:28)$) ($(90:2.75)+(0:6)+(\b:13)+(\c:28)$) -- ($(90+120*2:2.75)+(60:6)+(\b:13)+(\c:28)$);
\foreach \a in {0,60} {
\draw ($(90:2.75)+(\a:6)+(\b:13)+(\c:28)$) --  ($(90+120:2.75)+(\a:6)+(\b:13)+(\c:28)$) -- ($(90+120*2:2.75)+(\a:6)+(\b:13)+(\c:28)$) -- ($(90:2.75)+(\a:6)+(\b:13)+(\c:28)$) ($(90:2.75)+(120*2:2)+(\a:6)+(\b:13)+(\c:28)$) to ++ ($(0:2)$) ($(90+120:2.75)+(0:2)+(\a:6)+(\b:13)+(\c:28)$) to ++ ($(120:2)$)  ($(90+120*2:2.75)+(120:2)+(\a:6)+(\b:13)+(\c:28)$) to ++ ($(120*2:2)$) ;
}
}
\draw ($(90:2.75)+(60:6)+(60:13)$) -- ($(90+120:2.75)+(60:6)+(60:13)+(60:28)$) ($(90+120:2.75)+(0:6)+(0:13)$) -- ($(90+120:2.75)+(0:6)+(0:13)+(0:28)$) ($(90:2.75)+(60:6)+(0:13)+(0:28)$) -- ($(90+120*2:2.75)+(0:6)+(0:13)+(60:28)$);
}
\foreach \d in {0,60} {
\draw ($(90:2.75)+(\d:60)$) --  ($(90+120:2.75)+(\d:60)$) -- ($(90+120*2:2.75)+(\d:60)$) -- ($(90:2.75)+(\d:60)$) ($(90:2.75)+(120*2:2)+(\d:60)$) to ++ ($(0:2)$) ($(90+120:2.75)+(0:2)+(\d:60)$) to ++ ($(120:2)$)  ($(90+120*2:2.75)+(120:2)+(\d:60)$) to ++ ($(120*2:2)$) ;
\draw ($(90:2.75)+(\d:60)$) -- ($(90+120:2.75)+(60:6)+(\d:60)$) ($(90+120*2:2.75)+(\d:60)$) -- ($(90+120:2.75)+(0:6)+(\d:60)$) ($(90:2.75)+(0:6)+(\d:60)$) -- ($(90+120*2:2.75)+(60:6)+(\d:60)$);
\foreach \a in {0,60} {
\draw ($(90:2.75)+(\a:6)+(\d:60)$) --  ($(90+120:2.75)+(\a:6)+(\d:60)$) -- ($(90+120*2:2.75)+(\a:6)+(\d:60)$) -- ($(90:2.75)+(\a:6)+(\d:60)$) ($(90:2.75)+(120*2:2)+(\a:6)+(\d:60)$) to ++ ($(0:2)$) ($(90+120:2.75)+(0:2)+(\a:6)+(\d:60)$) to ++ ($(120:2)$)  ($(90+120*2:2.75)+(120:2)+(\a:6)+(\d:60)$) to ++ ($(120*2:2)$);
}
\draw ($(90:2.75)+(60:6)+(\d:60)$) -- ($(90+120:2.75)+(60:6)+(60:13)+(\d:60)$) ($(90+120*2:2.75)+(0:6)+(\d:60)$) -- ($(90+120:2.75)+(0:6)+(0:13)+(\d:60)$) ($(90:2.75)+(60:6)+(0:13)+(\d:60)$) -- ($(90+120*2:2.75)+(0:6)+(60:13)+(\d:60)$);
\foreach \b in {0,60} {
\draw ($(90:2.75)+(\b:13)+(\d:60)$) --  ($(90+120:2.75)+(\b:13)+(\d:60)$) -- ($(90+120*2:2.75)+(\b:13)+(\d:60)$) -- ($(90:2.75)+(\b:13)+(\d:60)$) ($(90:2.75)+(120*2:2)+(\b:13)+(\d:60)$) to ++ ($(0:2)$) ($(90+120:2.75)+(0:2)+(\b:13)+(\d:60)$) to ++ ($(120:2)$)  ($(90+120*2:2.75)+(120:2)+(\b:13)+(\d:60)$) to ++ ($(120*2:2)$) ;
\draw ($(90:2.75)+(\b:13)+(\d:60)$) -- ($(90+120:2.75)+(60:6)+(\b:13)+(\d:60)$) ($(90+120*2:2.75)+(\b:13)+(\d:60)$) -- ($(90+120:2.75)+(0:6)+(\b:13)+(\d:60)$) ($(90:2.75)+(0:6)+(\b:13)+(\d:60)$) -- ($(90+120*2:2.75)+(60:6)+(\b:13)+(\d:60)$);
\foreach \a in {0,60} {
\draw ($(90:2.75)+(\a:6)+(\b:13)+(\d:60)$) --  ($(90+120:2.75)+(\a:6)+(\b:13)+(\d:60)$) -- ($(90+120*2:2.75)+(\a:6)+(\b:13)+(\d:60)$) -- ($(90:2.75)+(\a:6)+(\b:13)+(\d:60)$) ($(90:2.75)+(120*2:2)+(\a:6)+(\b:13)+(\d:60)$) to ++ ($(0:2)$) ($(90+120:2.75)+(0:2)+(\a:6)+(\b:13)+(\d:60)$) to ++ ($(120:2)$)  ($(90+120*2:2.75)+(120:2)+(\a:6)+(\b:13)+(\d:60)$) to ++ ($(120*2:2)$) ;
}
}
\foreach \c in {0,60} {
\draw ($(90:2.75)+(\c:28)+(\d:60)$) --  ($(90+120:2.75)+(\c:28)+(\d:60)$) -- ($(90+120*2:2.75)+(\c:28)+(\d:60)$) -- ($(90:2.75)+(\c:28)+(\d:60)$) ($(90:2.75)+(120*2:2)+(\c:28)+(\d:60)$) to ++ ($(0:2)$) ($(90+120:2.75)+(0:2)+(\c:28)+(\d:60)$) to ++ ($(120:2)$)  ($(90+120*2:2.75)+(120:2)+(\c:28)+(\d:60)$) to ++ ($(120*2:2)$) ;
\draw ($(90:2.75)+(\c:28)+(\d:60)$) -- ($(90+120:2.75)+(60:6)+(\c:28)+(\d:60)$) ($(90+120*2:2.75)+(\c:28)+(\d:60)$) -- ($(90+120:2.75)+(0:6)+(\c:28)+(\d:60)$) ($(90:2.75)+(0:6)+(\c:28)+(\d:60)$) -- ($(90+120*2:2.75)+(60:6)+(\c:28)+(\d:60)$);
\foreach \a in {0,60} {
\draw ($(90:2.75)+(\a:6)+(\c:28)+(\d:60)$) --  ($(90+120:2.75)+(\a:6)+(\c:28)+(\d:60)$) -- ($(90+120*2:2.75)+(\a:6)+(\c:28)+(\d:60)$) -- ($(90:2.75)+(\a:6)+(\c:28)+(\d:60)$) ($(90:2.75)+(120*2:2)+(\a:6)+(\c:28)+(\d:60)$) to ++ ($(0:2)$) ($(90+120:2.75)+(0:2)+(\a:6)+(\c:28)+(\d:60)$) to ++ ($(120:2)$)  ($(90+120*2:2.75)+(120:2)+(\a:6)+(\c:28)+(\d:60)$) to ++ ($(120*2:2)$);
}
\draw ($(90:2.75)+(60:6)+(\c:28)+(\d:60)$) -- ($(90+120:2.75)+(60:6)+(60:13)+(\c:28)+(\d:60)$) ($(90+120*2:2.75)+(0:6)+(\c:28)+(\d:60)$) -- ($(90+120:2.75)+(0:6)+(0:13)+(\c:28)+(\d:60)$) ($(90:2.75)+(60:6)+(0:13)+(\c:28)+(\d:60)$) -- ($(90+120*2:2.75)+(0:6)+(60:13)+(\c:28)+(\d:60)$);
\foreach \b in {0,60} {
\draw ($(90:2.75)+(\b:13)+(\c:28)+(\d:60)$) --  ($(90+120:2.75)+(\b:13)+(\c:28)+(\d:60)$) -- ($(90+120*2:2.75)+(\b:13)+(\c:28)+(\d:60)$) -- ($(90:2.75)+(\b:13)+(\c:28)+(\d:60)$) ($(90:2.75)+(120*2:2)+(\b:13)+(\c:28)+(\d:60)$) to ++ ($(0:2)$) ($(90+120:2.75)+(0:2)+(\b:13)+(\c:28)+(\d:60)$) to ++ ($(120:2)$)  ($(90+120*2:2.75)+(120:2)+(\b:13)+(\c:28)+(\d:60)$) to ++ ($(120*2:2)$) ;
\draw ($(90:2.75)+(\b:13)+(\c:28)+(\d:60)$) -- ($(90+120:2.75)+(60:6)+(\b:13)+(\c:28)+(\d:60)$) ($(90+120*2:2.75)+(\b:13)+(\c:28)+(\d:60)$) -- ($(90+120:2.75)+(0:6)+(\b:13)+(\c:28)+(\d:60)$) ($(90:2.75)+(0:6)+(\b:13)+(\c:28)+(\d:60)$) -- ($(90+120*2:2.75)+(60:6)+(\b:13)+(\c:28)+(\d:60)$);
\foreach \a in {0,60} {
\draw ($(90:2.75)+(\a:6)+(\b:13)+(\c:28)+(\d:60)$) --  ($(90+120:2.75)+(\a:6)+(\b:13)+(\c:28)+(\d:60)$) -- ($(90+120*2:2.75)+(\a:6)+(\b:13)+(\c:28)+(\d:60)$) -- ($(90:2.75)+(\a:6)+(\b:13)+(\c:28)+(\d:60)$) ($(90:2.75)+(120*2:2)+(\a:6)+(\b:13)+(\c:28)+(\d:60)$) to ++ ($(0:2)$) ($(90+120:2.75)+(0:2)+(\a:6)+(\b:13)+(\c:28)+(\d:60)$) to ++ ($(120:2)$)  ($(90+120*2:2.75)+(120:2)+(\a:6)+(\b:13)+(\c:28)+(\d:60)$) to ++ ($(120*2:2)$) ;
}
}
\draw ($(90:2.75)+(60:6)+(60:13)+(\d:60)$) -- ($(90+120:2.75)+(60:6)+(60:13)+(60:28)+(\d:60)$) ($(90+120:2.75)+(0:6)+(0:13)+(\d:60)$) -- ($(90+120:2.75)+(0:6)+(0:13)+(0:28)+(\d:60)$) ($(90:2.75)+(60:6)+(0:13)+(0:28)+(\d:60)$) -- ($(90+120*2:2.75)+(0:6)+(0:13)+(60:28)+(\d:60)$);
}
\draw ($(90:2.75)+(60:6)+(60:13)+(60:28)$) -- ($(90+120:2.75)+(60:6)+(60:13)+(60:28)+(60:60)$) ($(90+120:2.75)+(0:6)+(0:13)+(0:28)$) -- ($(90+120:2.75)+(0:6)+(0:13)+(0:28)+(0:60)$) ($(90:2.75)+(60:6)+(0:13)+(0:28)+(0:60)$) -- ($(90+120*2:2.75)+(0:6)+(0:13)+(0:28)+(60:60)$);
}
\end{tikzpicture}
\caption{The Sierpinski gasket (SG) and the stretched Sierpinski gasket (SSG).}
\label{SG and SSG}
\end{figure}

Nevertheless, these models can fail to capture the essential structure of the media they are supposed to describe. The main message of the present paper is that reducing the analysis on an object to one-dimensional analysis on a union of lines can ignore a significant part of its intrinsic structure and therefore give a far too simple, hence incomplete, framework to investigate analytical questions on it. We aim to furnish the latter statement by studying here what we call the \textit{stretched Sierpinski gasket}, SSG for short, in $\BbR^2$. This space has also been called the ``deformed Sierpinski gasket''~\cite{OSS90} or ``Hanoi attractor''~\cite{ARF12, ARKT16, ARF16} and it is obtained from the classical Sierpinski gasket SG by replacing each branching point of the SG by an interval (see Figure~\ref{SG and SSG}).

\medskip

As a result, SSG is the closure of a countable union of one-dimensional intervals. One could thus think of constructing and analysing diffusion processes on it via quantum graphs/cable systems, an approach that has actually been considered in several works~\cite{GK14,ARKT16}. 

 
Let us give a rough definition of a cable system/quantum graph, leaving details to~\cite{BB04,Kuc08}. Starting from a weighted graph $(V, E, C)$ with vertex set $V$, edge set $E \subseteq \{(p, q)~|~p, q \in V\}$ and edge conductances/weights $C = \{C_{pq}~|~(p, q) \in E\}$, each edge $(p, q) \in E$ is identified with the line segment parametrized by $\xi_{p,q}(t) = (1 - t)p + tq$, $t \in [0, 1]$, and equipped with the Dirichlet energy $\D_{pq}$ on the line segment $pq$ given by
\[
\D_{pq}(\cdot, \cdot) = \int_0^1 \frac{d(\cdot\circ\xi_{p, q})}{dt}\frac{d(\cdot \circ\xi_{p, q})}{dt}dt.
\]
The consequent energy form $\E$ on the whole space is thus defined as
\[
\E(u, v) = \sum_{pq\in E} C_{pq}\D_{pq}(u, v),
\]
where the domain of $\E$ consists of all continuous $L^2$-functions on the whole space whose restriction to each edge $pq$ belongs to the Sobolev space $H^1(\xi_{p,q}([0,1]),dx)$. In a natural way, this quadratic form $\E$ induces a diffusion process that behaves like one-dimensional Brownian motion on each edge.

\medskip

Following this direction, a diffusion on SSG might be expected to consist basically in gluing the different Brownian motions on each interval. However, in considering SSG as a union of one-dimensional lines, one overlooks the ``reminiscence'' of SG in the geometric structure of SSG. In fact, the cable system/quantum graph approach disregards the underlying geometry of the space in the sense that it ignores the considerable role played by the arrangement of the vertices in space. Furthermore, classical quantum graph theory requires some finiteness condition that makes it inapplicable to cases such as fractals or infinite trees.

\medskip

Indeed, we show in this paper that the geometric ``reminiscence'' of the Sierpinski gasket also appears in the diffusion on SSG, a fact that stays hidden when using cable systems/quantum graphs. 

\medskip

The diffusion processes considered here will be associated with a Dirichlet form induced by a \textit{completely symmetric resistance form}. The theory of resistance forms was introduced in~\cite{Kig01} and further developed in particular to study analysis on ``low-dimensional'' fractals from an intrinsic point of view, see~\cite{Kig12} and references therein. Their most representative property is that, unlike Dirichlet forms, they are defined without requiring any measure on the underlying space. In our case, a completely symmetric resistance form $(\E,\F)$ on SSG is a resistance form whose restriction to every contractive copy in itself is invariant under all geometrical symmetries of the copy. More precisely, let $X$ be a subset of SSG which is similar to SSG itself and let $G: SSG \to X$ be the associated contractive similitude. If we denote by $\E_X$ the part of the original form $\E$ associated with $X$, then
\[
\E_X(u{\circ}G^{-1}, v{\circ}G^{-1}) 
\]
is again a form on SSG. We say that $(\E,\F)$ is completely symmetric if $\E_X(u\circ G^{-1}, v\circ G^{-1})$ is invariant under any isometry of the regular triangle. (See Section~\ref{RFF} for the exact definition.)

\medskip

As a key step towards the study of such diffusion processes, the present paper is devoted to establishing the existence of completely symmetric resistance forms on SSG. Even more, we provide a full characterization of all possible forms of this type by showing in Theorem~\ref{RFF.thm00} that any completely symmetric resistance form on SSG can be written as
\[
a\E^*(\cdot,\cdot)+b\D^I_\eta(\cdot,\cdot)
\]
for some $a \ge 0$ and $b > 0$. The forms $\E^*$ and $\D^I_{\eta}$ are briefly explained below. Conversely, we will show that any linear combination of $\E^*$ and $\D_\eta^I$ as above with $a \ge 0$ and $b > 0$ can be realized as a resistance form on SSG.

\medskip

On the one hand, $\D^I_{\eta}$ arises as a limit of sums of standard Dirichlet energies and it is defined as follows. Let $\eta = \{\eta_m\}_{m \ge 1} \subseteq (0, \infty)$ satisfy $\sum_{m \ge 1} (\frac 53)^{m - 1}\eta_m = 1$ and let $\D_k^I$ be the sum of the Dirichlet integrals over the line segments that appear in the $k$th approximation step of SSG for the first time, i.e.
\[
\D^I_k(u, v) = \sum_{pq \in J_k\setminus J_{k-1}} \D_{pq}(u, v),
\]
where $J_k$ denotes the set of line segments in the $k$-th approximation step. 
The quadratic form $\D^I_{\eta}$ is defined as the weighted sum of the $\D^I_k$'s whose weights are given by $\eta = \{\eta_m\}_{m \ge 1}$, i.e.
\[
\D^I_{\eta}(u, v) = \sum_{k\geq 1}\frac{1}{\eta_k}\D_k^I(u,v).
\]
It resembles the cable system/quantum graph approach in this setting. In particular, the special case $a = 0$ has been called ``fractal quantum graph'' in~\cite{ARKT16}, where the authors have shown that $\D^I_{\eta}$ is a resistance form for some limited choices of $\eta$. 

\medskip

On the other hand, the form $\E^*$ corresponds to the standard resistance form on SG. (See Definition~\ref{SG.def10} and~\cite{Kig01} for further details about this form.) Notice that any function on SG can be thought of as a function on SSG by making its value constant on each line segment. In this manner, we can regard the standard resistance form on SG as a quadratic form on SSG, see Definition~\ref{RFF.def02} for a precise formulation. This part of $\E$, which may be called the ``fractal part'', had remained unseen in the previous works~\cite{GK14, ARKT16} because there only limits of quantum graphs were considered.

\medskip

In conclusion, this paper reveals that SSG is more than just the combination of a countably infinite number of line segments, not only from a geometric, but also from an analytic point of view, since the reminiscence of the Sierpinski gasket in SSG remains essentially present in both of them.

\medskip

We will begin our exposition by discussing the geometry of SSG in Section~\ref{GK}, providing a detailed construction as well as some of its most relevant intrinsic geometric properties. Section~\ref{SG} reviews the construction of the standard resistance form on SG and establishes a first link between functions on SG and on SSG. Completely symmetric resistance forms on SSG are rigorously introduced in Section~\ref{RFF} and the main classification result of this paper is stated in Theorem~\ref{RFF.thm00}. The forthcoming sections develop the machinery to prove this theorem: Section~\ref{CRF} proceeds with the construction of resistance forms on SSG by means of compatible sequences based on sequences of what we call \textit{matching pairs} of resistances. We will see in Section~\ref{IRF} that any completely symmetric resistance form on SSG actually corresponds to a constant multiple of a resistance form on SSG derived from a sequence of matching pairs. Once this correspondence is settled, Section~\ref{CMP} establishes a preliminary classification result for resistance forms $(\E_\R,\F_\R)$ derived from matching pairs displayed in Theorem~\ref{CRF.thm10}. At this point, any such form $\E_\R$ becomes the sum of an \textit{SG part} and a \textit{line part}. In this way, the reminiscence of SG in SSG comes to light. In Section~\ref{SRF}, the previous theorem is enhanced through a projection mapping onto the resistance forms having only line part. Section~\ref{DRF} is the core of the paper: in Theorem~\ref{DRF.thm10}, the domain of the completely symmetric resistance forms on SSG is fully described, and the SG part and the line part get their corresponding expression as the aforementioned forms $\E^*$, $\D_\eta^I$ respectively. This characterization will finally lead to the classification provided by Theorem~\ref{RFF.thm00}.

\section{Glossary of notations}

For the convenience of the reader we give below an index which summarizes notation repeatedly used throughout the text, and where the definitions may be found.

\begin{multicols}{2}
\noindent$B = \{(1, 2), (2, 3), (3, 1)\}$\\
$d_E$: the restriction of the Euclidean metric\\
$\D_{pq}$ -- Definition~\ref{RFF.def00}\\
$\D_m^I$ -- Definition~\ref{RFF.def00}\\
$\D_{\eta}^I$ -- Definition~\ref{RFF.def15}\\
$e_{ij}$ -- Definition~\ref{GK.def10}\\
$e_{ij}^w = G_w(e_{ij})$\\
$E_{\R, m}$ -- Definition~\ref{RNS.def20}\\
$\E_m^*$ -- Definition~\ref{SG.def10}\\
$\E_*$ -- Theorem~\ref{SG.thm10}\\
$\E^*$ -- Definition~\ref{RFF.def02}\\
$\E|_Y$ -- Proposition~\ref{BRF.prop20}\\
$\E_{\R}$ -- Definition~\ref{RNS.def10}\\
$\E_{\R, m}$ -- Definition~\ref{RNS.def500}\\
$\E_{\R}^I$ -- \eqref{IND.eq100}\\
$\E_{\R}^{\SS}$ -- \eqref{IND.eq200}\\
$\widehat{\E}_{\R}$ -- Definition~\ref{RNS.def20}\\
$F_i = G_i$ in case of $\a = 0$\\
$\F^*$ -- Theorem~\ref{SG.thm10}\\
$\F_{\R}$ -- Definition~\ref{RNS.def10}\\
$\F_{\R}^I$ -- Definition~\ref{CMP.def10}\\
$\F_{\R}^{\SS}$ -- Definition~\ref{CMP.def10}\\
$\F_{\R, *}$ -- Definition~\ref{CRF.def00}\\
$\F^{\SS}$ -- Definition~\ref{RFF.def02}\\
$\F^*_m$ -- Definition~\ref{RFF.def16}\\
$\F_{\infty}^*$ -- Definition~\ref{RFF.def16}\\
$\F_{\eta}$ -- Definition~\ref{RFF.def15}\\
$\F_\eta^*$ -- Definition~\ref{RFF.def15}\\
$\F^{(n)}$ -- Definition~\ref{DRF.def10}\\
$\F|_Y$ -- Definition~\ref{BRF.def30}\\
$\F_0(Y)$ -- Definition~\ref{BRF.def30}\\
$\tF_m$ -- Definition~\ref{RFF.def16}\\
$\tF_{\infty}$ -- Definition~\ref{RFF.def16}\\
$\widetilde{\F}$ -- Definition~\ref{RFF.def00}\\
$\widehat{\F}_{\R}$ -- Definition~\ref{RNS.def20}\\
$G_i$ -- Definition~\ref{GK.def10}\\
$G_w$ -- Definition~\ref{GK.def20}\\
$\G_K$ -- \eqref{IND.GK}\\
$h_Y$ -- Proposition~\ref{BRF.prop20}\\
$H^1([0, 1])$ -- Definition~\ref{RFF.def00}\\
$H^1(pq)$ -- Definition~\ref{RFF.def00}\\
$\H_{(\E, \F)}(Y)$ -- Definition~\ref{BRF.def40}\\
$K = K_{\a}$ for $\a \in (0, 1)$\\
$K_* = K_0 = $\,\,the Sierpinski gasket\\
$K_{\a}$ -- Proposition~\ref{GK.prop10}\\
$\L(\R)$ -- Definition~\ref{SRF.def10}\\
$\ell(V) = \{u| u: V \to \BbR\}$\\
$\MP$ -- Definition~\ref{CRF.def30}\\
$(\MP^{\BbN})^I$ -- \eqref{IND.eq300}\\
$\{p_1, p_2, p_3\}$: vertices of a regular triangle\\
$p_{ij}$ -- Definition~\ref{GK.def10}\\
$Q_0$ -- Definition~\ref{SG.def10}\\
$Q_m^I$ -- Definition~\ref{CRF.def20}\\
$Q_m^{\SS}$ -- Definition~\ref{RFF.def02}\\
$r_0(\E, \F)$ -- Definition~\ref{IRF.def10}\\
$R_*$ -- Theorem~\ref{DRF.thm10}\\
$R_*^{(n)}$ -- Definition~\ref{DRF.def10}\\
$\R_{(\E, \F)}$ -- Definition~\ref{IRF.def20}\\
$\R^{(n)}$ -- Definition~\ref{DRF.def10}\\
$\RF_S$ -- Definition~\ref{RFF.def10}\\
$\RF^{(0)}_S$ -- Definition~\ref{RFF.def10}\\
$\RF_S^N$ -- Definition~\ref{IRF.def10}\\
$S = \{1, 2, 3\}$\\
$V_m$ -- Definition~\ref{GK.def20}\\
$V_m^*$ -- Definition~\ref{SG.def10}\\
$W_m$ -- Definition~\ref{GK.def20}\\
$W_*$ -- Definition~\ref{GK.def20}\\
$\alpha$ -- Definition~\ref{GK.def10}\\
$\iota = \iota^{\a}$\\
$\iota^{\a}$ -- Proof of Proposition~\ref{GK.prop20}\\
$\pi, \pi_*, \pi^*$ -- Definition~\ref{def05}\\
$\eta^{(n)}$ -- Definition~\ref{DRF.def10}\\
$\eta^{(n)}_m$ -- Definition~\ref{DRF.def10}\\
$\SS = \SS_{\a}$ for $\a \in (0, 1)$\\
$\SS_{\a}$ -- Definition~\ref{GK.prop10}\\
$\xi_{pq}$ -- Definition~\ref{RFF.def00}
\end{multicols}
\setcounter{equation}{0}

\section{Geometry of $K$}\label{GK}
In this section, we set up the geometric construction of SSG in $\BbR^2$ and fix the corresponding notation that will be carried throughout the paper. \par
Let $S = \{1, 2, 3\}$ and let $\{p_1, p_2, p_3\}$ be the collection of vertices of a regular triangle in $\BbR^2$. For the purpose of normalization, we assume that $p_1 + p_2 + p_3 = 0$ and $|p_i - p_j| = 1$ for any $i \neq j$. 

\begin{figure}[H]
\centering
\begin{tikzpicture}[scale=0.4]
\coordinate[label=below:{$p_2$}] (p_1) at (0,0);
\coordinate[label=above:{$p_1$}] (p_2) at (3.5,6.0621);
\coordinate[label=below:{$p_3$}] (p_3) at (7,0);
\coordinate (p_12) at (1.5,2.59805);
\coordinate (p_21) at (2,3.46405);
\coordinate (p_13) at (3,0);
\coordinate (p_31) at (4,0);
\coordinate (p_23) at (5,3.46405);
\coordinate (p_32) at (5.5,2.59805);
\node at (3.5,-1.5) {$V_0$};

\draw (p_1) -- (p_2) node[midway, left] {$1$} -- (p_3) node[midway, right] {$1$} -- (p_1) node[midway, above] {$1$};
\fill (p_1) circle (5pt);
\fill (p_2) circle (5pt);
\fill (p_3) circle (5pt);
\end{tikzpicture}
\hspace*{.5in} 
\begin{tikzpicture}[scale=0.4]
\coordinate (p_1) at (0,0);
\coordinate (p_2) at (3.5,6.0621);
\coordinate (p_3) at (7,0);
\coordinate[label=left:{$p_{21}$}] (p_12) at (0.875,1.5155);
\coordinate[label=left:{$p_{12}$}] (p_21) at (2.625,4.5466);
\coordinate[label=below:{$p_{23}$}] (p_13) at (1.75,0);
\coordinate[label=below:{$p_{32}$}] (p_31) at (5.25,0);
\coordinate[label=right:{$p_{13}$}] (p_23) at (4.375,4.5466);
\coordinate[label=right:{$p_{31}$}] (p_32) at (6.125,1.5155);

\node at (3.5,-1.5) {$V_1$};

\draw (p_1) -- (p_2) -- (p_3) -- (p_1) (p_13) -- (p_12) (p_21) -- (p_23) (p_31) -- (p_32);
\fill (p_1) circle (5pt);
\fill (p_2) circle (5pt);
\fill (p_3) circle (5pt);
\fill (p_12) circle (5pt);
\fill (p_21) circle (5pt);
\fill (p_31) circle (5pt);
\fill (p_13) circle (5pt);
\fill (p_23) circle (5pt);
\fill (p_32) circle (5pt);
\end{tikzpicture}
\caption{Geometric construction}
\label{Fig1}
\end{figure}
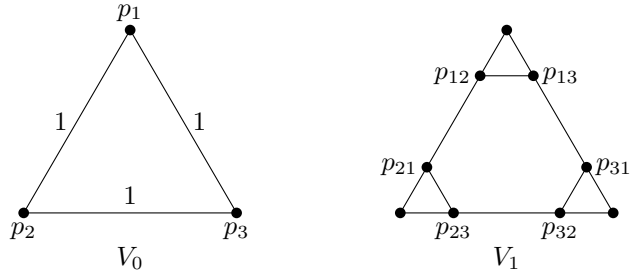

\definition\label{GK.def10}
For each $i \in S$, define $G_i\colon\BbR^2\to\BbR^2$ by
\[
G_i(x) = \frac{1 - \a}2(x - p_i) + p_i,
\]
where $0 \le \a \le 1$. Moreover, set $p_{ij} = G_i(p_j)$ for $i \neq j$ and denote by $e_{ij}$ the line segment $p_{ij}p_{ji}$.
\enddefinition

If $\a = 0$, then $p_{ij} = p_{ji}$ for any $i \neq j$ and hence $e_{ij} = \{p_{ij}\}$. Notice that $G_i$, $p_{ij}$ and $e_{ij}$ actually depend on $\a$. However, we will see in Proposition~\ref{GK.prop20} that the sets $K_{\a}$ defined in Propostion~\ref{GK.prop10} are homeomorphic to each other for $\a \in (0, 1)$ and therefore we do not write $\a$ explicitly in the notation.

\definition\label{GK.def20}
Let $W_0 = \{\emptyset\}$ and define
\[
W_m = S^m = \{w~|~w = w_1\ldots{w_m}, w_i \in S~\text{for any } i = 1, \ldots, m\}
\]
for $m \ge 1$, as well as $W_* = \cup_{m \ge 0} W_m$. Moreover, for any $w = w_1\ldots{w_m} \in W_*$, define $G_w\colon\BbR^2\to\BbR^2$ by
\[
G_w = G_{w_1}{\circ}G_{w_2}\circ\cdots{\circ}G_{w_m}
\]
and $G_\emptyset=\operatorname{id}$ is the identity map on $\BbR^2$. 
Finally, set $V_0 = \{p_1, p_2, p_3\}$ and
\[
V_m = \bigcup_{w \in W_m} G_w(V_0)
\]
for $m \ge 1$.
\enddefinition

\notation
From now on, we denote by $B = \{(1, 2), (2, 3), (3, 1)\}$, where $B$ stands for the word ``Bond'', and write $e_{ij}^w = G_w(e_{ij})$ for any $(w ,(i, j)) \in W_*\times B$.
\endnotation

\prop\label{GK.prop10}
For any $0\leq \a\leq 1$, there exists a unique compact set $K_{\a}\subseteq\BbR^2$ such that
\[
K_{\a} = G_1(K_{\a}) \cup G_2(K_{\a}) \cup G_3(K_{\a}) \cup e_{12} \cup e_{23} \cup e_{31}.
\]
Furthermore,
\[
K_{\a} = \SS_{\a} \cup \bigcup_{(w,(i,j)) \in W_*\times B} e_{ij}^w,
\]
where $\SS_{\a}$ is the self-similar set associated with $\{G_1, G_2, G_3\}$, i.e. $\SS_{\a}$ is the unique nonempty compact set satisfying
\begin{equation}\label{GK.eq100}
\SS_{\a} = G_1(\SS_{\a}) \cup G_2(\SS_{\a}) \cup G_3(\SS_{\a}).
\end{equation}
Moreover, $\cup_{m \ge 0} V_m$ is a dense subset of $\SS_{\a}$.
\endprop
\begin{proof}
This follows from~\cite[Section 4, Theorem 1]{Hat85} since $G_1,G_2,G_3$ are $\frac{1-\alpha}{2}$-contractions.
\end{proof}

\remark
$\SS_{\a}$ is a Cantor set for any $0 < \a < 1$.
\endremark

Notice that $K_0$ coincides with the Sierpinski gasket while $K_1$ is the union of the three line segments $p_1p_2$, $p_2p_3$ and $p_3p_1$. Whenever $\a \in (0,1)$, we can refer to any of $K_\a$ as \textit{the} stretched Sierpinski gasket SSG in view of the next proposition.
\prop\label{GK.prop20}
The sets $K_\a$, $\a \in (0,1)$, are pairwise homeomorphic.
\endprop

\begin{proof}
Use $G_w^{\a}$ and $e_{ij}^{\a, w}$ to denote $G_w$ and $e_{ij}^w$ respectively. Note that $\SS_{\a}$ is naturally homeomorphic to $\{1, 2, 3\}^{\BbN}$ by the canonical coding map $\iota^{\a}$ defined by $\{\iota^{\a}(\omega_1\omega_2\ldots)\} = \cap_{m \ge 1} G^w_{\omega_1\ldots\omega_m}(K_{\a})$. Let $\vp_{\a_1, \a_2} = \iota^{\a_2}{\circ}(\iota^{\a_1})^{-1}$. Then $\vp_{\a_1, \a_2}: \SS_{\a_1} \to \SS_{\a_2}$ is a homeomorphism. Extend $\vp_{\a_1, \a_2}$ onto $e_{ij}^{w, \a_1}$ by $\vp_{\a_1, \a_2}|_{e_{ij}^{w, \a_1}} = G_w^{\a_2}\circ(G_w^{\a_1})^{-1}|_{e_{ij}^{w, \a_1}}$ for any $i, j \in B$ and $w \in W_*$. Then $\vp_{\a_1, \a_2}: K_{\a_1} \to K_{\a_2}$ is a homeomorphism.
\end{proof}

Since resistance forms on $K_\a$ only depend on the topological structure of $K_{\a}$, which is the same for any $\a \in (0, 1)$ due to the previous proposition, we will omit $\a$ in the definition given by Proposition~\ref{GK.prop10} and write $K = K_{\a}$ and $\SS = \SS_{\a}$ as long as $\a \in (0, 1)$. Moreover, we will consider $d_E$ to be the restriction of the Euclidean metric to $K_{1/2}$ and regard $d_E$ as the canonical metric on $K$.

In view of~\eqref{GK.eq100}, there exists a canonical map $\iota\colon S^{\BbN} \to \SS$ defined by $\iota(\omega_1\omega_2\ldots) = \bigcap_{m \ge 1} G_{\word{\omega}m}(\SS)$. Through this map $\iota$, we identify $\SS$ with $S^{\BbN}$ hereafter in this paper.

\setcounter{equation}{0}

\section{The Sierpinski gasket}\label{SG}

As already mentioned, if $\a = 0$ in Definition~\ref{GK.def10}, then $K_{\a}$ is the Sierpinski gasket, $p_{ij} = p_{ji}$ and $e_{ij} = \{p_{ij}\}$ for any $(i, j) \in B$. In this case, we will denote $G_i$ and $K_\a$ by $F_i$ and $K_*$ respectively. We explain in this section how to view continuous functions on the Sierpinski gasket $K_*$ as continuous functions on the stretched Sierpinski gasket $K$ and review the construction of the standard resistance form on $K_*$. Further details and proofs can be found e.g. in~\cite{Kig01}.

\definition\label{def05}
Let $\pi\colon \SS \to K_*$ be the canonical coding map given by $\{\pi(\omega_1\omega_2\ldots) \}= \bigcap\limits_{m \ge 0} F_{\omega_1\ldots\omega_m}(K_*)$ and define $\pi_*\colon K \to K_*$ by
\[
\pi_*|_{\SS} = \pi
\]
and
\[
\pi_*(e_{ij}^w) = \pi(wi(j)^{\infty}) = \pi(wj(i)^{\infty})
\]
for any $(w,(i, j)) \in W_*\times B$. Furthermore, define $\pi^*\colon C(K_*) \to C(K)$ by $\pi^*(u) = u\circ{\pi_*}$.
\enddefinition
From this definition it follows that $u\in \pi^*(C(K_*))$ if and only if $u \in C(K)$ and $u|_{e_{ij}^w}$ is constant for each $(w,(i, j)) \in W_*\times B$, a fact stated in the next proposition. Moreover, $\pi^*$ is injective and it preserves the supremum norm. We will thus identify $C(K_*)$ with $\pi^*(C(K_*))$ and think of $C(K_*)$ as a subset of $C(K)$ in this manner. Thus we have the following proposition.

\prop\label{SG.prop05}
\[
C(K_*) = \{u~|~u \in C(K),~u|_{e_{ij}^w} \text{ is constant for any }(w, (i, j)) \in W_* \times B\}.
\]
\endprop

We finish this paragraph with some classical definitions and results concerning the standard resistance form on $K_*$ that will become relevant to state our main theorem. 
\notation
For any set $V$ we use the standard notation $\ell(V)=\{u~|~u\colon V\to\BbR\}$.
\endnotation
\definition\label{SG.def10}
Let $V_0^* = \{p_1, p_2, p_3\}$ and define $V_m^*$ inductively by $V_{m + 1}^* = \cup_{i = 1}^3 F_i(V_m^*)$ for $m \ge 0$. 
Furthermore, let the quadratic form $\E^*_m(\cdot, \cdot)$ on $\l(V_m^*)$ be defined as
\[
Q_0(u, u) = \sum_{(i, j) \in B} (u(p_i) - u(p_j))^2
\]
for $m=0$, and
\[
\E_m^*(u, u) = \Big(\frac 53\Big)^m\sum_{w \in W_m} Q_0(u{\circ}F_w, u{\circ}F_w)
\]
for $m \ge 1$.
\enddefinition

\prop\label{SG.prop10}
For any $u : K_* \to \BbR$ and any $m \ge 0$,
\[
\E_m^*(u|_{V_m^*}, u|_{V_m^*}) \le \E_{m+1}^*(u|_{V_{m + 1}^*}, u|_{V_{m + 1}^*})
\]
and $\lim\limits_{m \to \infty} \E_m^*(u|_{V_m^*}, u|_{V_m^*}) = 0$ if and only if $u$ is constant on $K_*$.
\endprop
\begin{proof}
This follows directly from Definition~\ref{SG.def10}.
\end{proof}

\thm\label{SG.thm10}
Define
\[
\F^* = \{u~|~u \in C(K_*), \lim_{m \to \infty}\E_m^*(u|_{V_m^*}, u|_{V_m^*}) < +\infty\}
\]
and
\[
\E_*(u, u) = \lim_{m \to \infty} \E_m^*(u|_{V_m^*}, u|_{V_m^*})
\]
for $u \in \F^*$. Then $\F^* \subseteq C(K_*)$ and $(\E_*, \F^*)$ is a resistance form on $K_*$.
\endthm
\begin{proof}
See Theorem~\ref{BRF.thm10} or~\cite[Theorem 3.13]{Kig12}.
\end{proof}

Analogously to $C(K_*)$, we will identify $\F^*$ with $\pi^*(\F^*)$ and thus regard $\F^*$ as a subset of $C(K)$.

\setcounter{equation}{0}

\section{Completely symmetric resistance forms}\label{RFF}

This section is devoted to giving a rigorous definition of completely symmetric resistance forms on SSG and presenting the main theorem of this paper, Theorem~\ref{RFF.thm00}, which provides their complete characterization and classification by means of the forms $\E^*$ and $\D_\eta^I$. The proof of Theorem~\ref{RFF.thm00} will require a suitable combination of the results obtained in the succeeding sections and it will therefore be presented at the end of Section~\ref{DRF}. We start by introducing some auxiliary notation and definitions. Recall that we write $K = K_{\a}$ for any $\a\in (0,1)$.

\definition\label{RFF.def00}
(1)~Let $H^1([0, 1])$ denote the Sobolev space
\begin{multline*}
H^1([0, 1])  = \Big\{u~\Big|~u\colon [0, 1] \to \BbR,~\frac{du}{dx} \in L^2([0, 1], dx), \\\text{where }\frac{du}{dx} \text{ is the derivative of }u \text{ in the sense of distributions}\Big\}.
\end{multline*}
(2)~For any $p, q \in \BbR^2$, let $pq$ denote the line segment with extreme points $p$ and $q$, and let $\xi_{p,q} \colon [0, 1] \to pq$ be given by $\xi_{p, q}(t) = (1 - t)p + tq$. We define 
\[
 H^1(pq) = \{u~|~u\colon pq \to \BbR,~u{\circ}\xi_{p, q} \in H^1([0, 1])\}
\]
and
\[
\D_{pq}(u, v) = \int_0^1 \frac{d(u{\circ}\xi_{p, q})}{dx}\frac{d(v{\circ}\xi_{p, q})}{dx}dx
\]
for any $u, v \in H^1(pq)$.

\noindent
(3)~Define
\[
\widetilde{\F} = \{u~|~u \in C(K)\text{ and }u|_{e^w_{ij}} \in H^1(e_{ij}^w)\text{ for any }(w, (i, j)) \in W_* \times B\}
\]
as well as 
\[
\D_m^I(u, v) = \sum_{(w, (i, j)) \in W_{m - 1} \times B} \D_{e_{ij}^w}(u|_{e_{ij}^w}, v|_{e_{ij}^w})
\]
for any $u, v \in \widetilde{\F}$ and $m\geq 1$.
\enddefinition

We introduce now the family of completely symmetric resistance forms on $K$ that play the central role in the classification theorem. 

Basic definitions and notation concerning resistance forms are reviewed in Section~\ref{BRF}, see also~\cite{Kig12}. First of all, consider the set of all linear mappings under which $K$ is invariant, i.e.
\begin{equation}\label{IND.GK}
\G_K = \{\vp~|~\vp\colon \BbR^2 \to \BbR^2~\text{ linear and such that } \vp(K) = K\}.
\end{equation}
Notice that this is in fact the dihedral group of symmetries of the triangle.

\definition\label{RFF.def10}
(1)~Let $\RF^{(0)}_S$ be the collection of resistance forms $(\E, \F)$ on $K$ satisfying the following three conditions (a), (b) and (c):
\begin{itemize}[leftmargin=.25in]
\item[(a)] $\F \subseteq C(K)$ and $u \in \F$ if and only if $u|_{G_i(K)} \in \F|_{G_i(K)}=\{v|_{G_i(K)}\,|\,v\in\F\}$ for any $i\in S$, and $u|_{e_{ij}} \in H^1(e_{ij})$ for any $(i, j) \in B$.
\item[(b)] Let $R$ be the resistance metric associated with $(\E, \F)$. Then the identity map from $(K, d_E)$ to $(K, R)$ is a homeomorphism.
\item[(c)] For any $\vp \in \G_K$ and $u \in \F$, $u{\circ}\vp \in \F$ and
\[
\E(u{\circ}\vp, u{\circ}\vp) = \E(u, u).
\]
\end{itemize}
(2)~Define $\RF_S$ to be the collection of resistance forms $(\E, \F) \in \RF^{(0)}_S$ with the following property:

\noindent
There exist a sequence $\{(\E_m, \F_m)\}_{m \ge 0} \subseteq \RF^{(0)}_S$ and a sequence $\{\eta_m\}_{m \ge 1} \subseteq (0, \infty)$ such that $(\E_0, \F_0) = (\E, \F)$, $ \F_m=\{u{\circ}G_i~|~u \in \F_{m - 1}\} $ for any $m \ge 1$ and $i \in S$, and
\begin{equation}\label{RFF.eq10}
\E_{m - 1}(u, v) = \sum_{i = 1}^3 \E_{m}(u{\circ}G_i, v{\circ}G_i) + \frac 1{\eta_m}\D^I_1(u, v)
\end{equation}
for any $m \ge 1$ and $u, v \in \F_{m - 1}$. The sequence $\{(\E_m, \F_m, \eta_m)\}_{m \ge 0} \subseteq \RF_S^{(0)} \times (0, \infty)$ is called the \textit{resolution of} $(\E, \F)$.
\enddefinition

\remark
Although $\eta_0$ is not needed in the previous definition, we will always set $\eta_0 = 1$ for the sake of formality.
\endremark

Applying \eqref{RFF.eq10} repeatedly, one immediately obtains the following proposition.

\prop\label{RFF.prop10}
Let $(\E, \F) \in \RF_S$. If $\{(\E_m, \F_m, \eta_m)\}_{m \ge 0} \subseteq \RF_S^{(0)} \times (0, \infty)$ is the resolution of $(\E, \F)$, then 
\begin{multline}\label{RFF.eq15}
\F = \{u~|~u{\circ}G_w \in \F_m~\text{for any }w \in W_m\text{ and }\\
 u|_{e_{ij}^w} \in H^1(e_{ij}^w)~\text{for any }(w, (i, j)) \in \big(\cup_{k = 0}^{m - 1}W_k\big) \times B\}.
\end{multline}
Moreover, for any $m \ge 1$ and any $u, v \in \F$,
\begin{equation}\label{RFF.eq20}
\E(u, v) = \sum_{w \in W_m}\E_m(u{\circ}G_w, v{\circ}G_w) + \sum_{k = 1}^{m} \frac 1{\eta_k}\D_k^I(u, v).
\end{equation}
\endprop

The next quadratic form resembles the classical resistance form $(\E_*,\F^*)$ on $K_*$ of Theorem~\ref{SG.thm10} and it will be precisely the ``fractal part'' missed by the cable system/quantum graph approach discussed in the introduction.
\definition\label{RFF.def02}
Let the quadratic form $Q_m^{\SS}$ on $\l(V_m)$ be given by
\[
Q_0^{\SS}(u, u) = \sum_{(i, j) \in B} (u(p_i) - u(p_j))^2
\]
for any $u\in\ell(V_0)$, and by
\[
Q_m^{\SS}(u, u) = \sum_{w \in W_m} Q_0(u{\circ}G_w, u{\circ}G_w)
\]
for $m \ge 1$ and any $u \in \l(V_m)$. Moreover, define
\[
\F^{\SS} = \bigg\{u~\bigg|~u\in C(K)\text{ and }\Big\{\Big(\frac 53\Big)^mQ_m^{\SS}(u, u)\Big\}_{m \ge 0}\text{ is a Cauchy sequence}\bigg\},
\]
as well as
\[
\E^*(u, u) =  \lim_{m \to \infty} \left(\frac 53\right)^mQ_m^{\SS}(u, u)
\]
for $u \in \F^{\SS}$.
\enddefinition

\definition\label{RFF.def16}
(1)~For any $m \ge 1$, let
\[
\tF_m = \{u~|~u \in \tF,~u|_{G_w(K)}\text{ is constant for any }w \in W_m\},
\] 
and
\[
\tF_{\infty} = \cup_{m \ge 1} \tF_m.
\]
(2)~For any $m \ge 1$, let
\[
\F^*_m = \{u~|~u \in \tF,~u{\circ}G_w \in \F^*\text{ for any }w \in W_m\}
\]
and
\[
\F_{\infty}^* = \cup_{m \ge 1} \F_m^*.
\]
\enddefinition

\remark
Notice that $\tF_m \subseteq \F_m^* \subseteq \F^{\SS}$ and $\F^* \subseteq C(K_*) \subseteq \tF$.
\endremark

Finally, we introduce the quadratic form $\D^I_{\eta}$ as the weighted sum of Dirichlet integrals whose weights are given by sequences $\{\eta_m\}_{m \ge 1}\subseteq (0,\infty)$. This form is the part that mirrors the cable system/quantum graph approach of an energy form on SSG. 

\definition\label{RFF.def15}
Let $\eta = \{\eta_m\}_{m \ge 1}$ be a sequence of positive numbers and for any $u \in \tF$, let
\[
\D_{\eta}^I(u, u) = \sum_{m=1}^\infty \frac1{\eta_m}\D^I_m(u, u).
\]
(Note that $\D_{\eta}^I(u, u)$ is well-defined if we allow the value $\infty$.) Moreover, define
\begin{multline*}
\F_{\eta} = \{u~|~u \in \tF,~\D_{\eta}^I(u, u) < +\infty\text{ and there exists }\{u_n\}_{n \ge 1} \subseteq \tF_{\infty}\text{ such that}\\
\lim_{n \to \infty}\D_{\eta}^I(u - u_n, u - u_n) = 0\text{ and }\lim_{n \to \infty}u_n(x) = u(x)\text{ for any }x \in K\},
\end{multline*}
as well as
\begin{multline*}
\F_\eta^* = \{u~|~u \in \tF \cap \F^{\SS},~\D_{\eta}^I(u, u) < +\infty\text{ and there exists }\{u_n\}_{n \ge 1} \subseteq \F_{\infty}^*\\
\text{such that }\lim_{n \to \infty}\E^*(u - u_n, u - u_n) = \lim_{n \to \infty} \D_{\eta}^I(u - u_n, u - u_n) = 0\\
\text{ and }\lim_{n \to \infty}u_n(x) = u(x)\text{ for any }x \in K\}.
\end{multline*}
\enddefinition

Our main result fully characterizes and identifies all resistance forms in $\RF_S$ by showing the correspondence between resistance forms on SSG that belong to $\RF_S$ and linear combinations of the forms $\E^*$ and $\D^I_\eta$.
\thm\label{RFF.thm00}
{\rm (1)}~$(\E, \F) \in \RF_S$ if and only if there exist $a \ge 0, b > 0$ and a sequence $\eta = \{\eta_m\}_{m \ge 1} \subseteq (0, \infty)$ such that
\begin{equation}\label{RFF.eq30}
 \sum_{m = 1}^{\infty}\Big(\frac 53\Big)^{m - 1}\eta_m = 1,
 \end{equation}
\begin{equation}\label{RFF.eq31}
\F = \begin{cases}
\F_{\eta}\quad&\text{if }a = 0,\\
\F_{\eta}^*\quad&\text{if }a > 0,
\end{cases}
\end{equation}
and
\[
\E(u, v) = a\E^*(u, v) + b\D_{\eta}^I(u, v)
\]
for any $u, v \in \F$.

\noindent
{\rm (2)}~If $\eta = \{\eta_m\}_{m \ge 1} \subseteq (0, \infty)$ satisfies \eqref{RFF.eq30}, then $\F_{\eta} \subseteq \F_{\eta}^*$ and
\[
\F_{\eta} = \{u~|~u \in \F_{\eta}^*,~\E^*(u, u) = 0\}.
\]
\endthm

\remark
As we mentioned in the introduction, the case $a=0$ was treated in~\cite{ARKT16} for a restricted type of sequence $\eta$. We would also like to emphasize that, even though at first sight one might want to apply the abstract result in that paper~\cite[Theorem 8.1]{ARKT16} in order to obtain (part of) Theorem~\ref{CRF.thm10}, it is not possible to do so in this setting since in particular the resistance metric associated with $(\E,\F)$ does not lead to a geodesic metric on SSG.
\endremark

\setcounter{equation}{0}

\section{Basics on resistance forms}\label{BRF}

For convenience of the reader, we give in this section a summary of definitions and basic facts from the theory of resistance forms used within the paper. A detailed and more extensive exposition of this theory can be found e.g. in~\cite{Kig01,Kig12}. 

\definition\label{BRF.def10}
Let $X$ be a set. A pair $(\E,\F)$ is called a resistance form on $X$ if it satisfies the following conditions (RF1) through
(RF5): \\
(RF1)~$\F$ is a linear subspace of $\l(X)=\{u\,|\,u: X \to \BbR\}$ containing constants and $\E$ is a non-negative symmetric quadratic form on
$\F$. $\E(u, u) = 0$ if and only if $u$ is constant on $X$.\\
(RF2)~Let $\sim$ be the equivalence relation on $\F$ defined by $u\sim v$ if and only if $u - v$ is constant on $X$. Then,
$(\F/\!\!\sim,\E)$ is a Hilbert space. \\
(RF3)~If $x \neq y$, then there exists $u \in \F$ such that $u(x) \neq u(y)$.\\
(RF4)~For any $p, q \in X$, 
\[
\sup\Big\{\frac{|u(p) - u(q)|^2}{\E(u,u)}~\big|~u \in \F, \E(u, u) > 0\Big\} 
\]
is finite. The above supremum is denoted by $R_{(\E, \F)}(p, q)$ and it is called the resistance metric on $X$ associated with the resistance form $(\E, \F)$.\\
(RF5)~For any $u\in\F$, $\overline{u} \in \F$ and $\E(\bar u,\bar u) \le\E(u,u)$, where $\overline{u}$ is defined by
\[
\overline{u}(p) = \begin{cases} 1 &\text{if }u(p) \ge 1,\\
                     u(p) &\text{if }0 < u(p) < 1,\\
                     0 &\text{if }u(p) \le 0. \end{cases}
\]
\enddefinition

\prop\label{BRF.prop10}
If $(\E, \F)$ is a resistance form on a set $X$, then the associated resistance metric $R_{(\E, \F)}(\cdot, \cdot)$ is a distance on $X$.
\endprop

If the set $X$ is finite, any resistance form on $X$ is a non-negative quadratic form on $\l(X)\times\l(X)$ that satisfies several conditions stated in the following lemma.

\lemma\label{BRF.lemma10}
Let $V$ be a finite set. Then $(\E, \l(V))$ is a resistance form on $V$ if and only if there exists $(C_{pq})_{p, q \in V}$ such that for any $p \neq q \in V$, $C_{pq} = C_{qp} \ge 0$ and there exist $m \ge 0$ and $(p_0, p_1, \ldots, p_m) \in V^{m + 1}$ such that $p_0 = p, p_m = q$ and $C_{p_ip_{i + 1}} > 0$ for any $i = 0, \ldots, m - 1$ and
\[
\E(u, v) = \frac 12\sum_{p, q \in V} C_{pq}(u(p) - u(q))(v(p) - v(q))
\]
for any $u \in \l(V)$.
\endlemma

If the set $X$ is infinite, in many cases a resistance form on $X$ is constructed by means of a suitable sequence of resistance forms on finite sets that approximate $X$ as Theorem~\ref{BRF.thm10} indicates.

\definition\label{BRF.def20}
Let $V$ and $U$ be finite sets satisfying $V \subseteq U$ and let $(\E_V, \l(V))$ and $(\E_U, \l(U))$ be resistance forms on $V$ and $U$ respectively. We write $(\E_V, \l(V)) \le (\E_U, \l(U))$ if and only if
\[
\E_V(u, u) = \min\{\E_U(v, v)~|~v \in \l(U),~v|_V = u\}
\]
for any $u \in \l(V)$. Let $V_m$ be a finite set and let $(\E_m, \l(V_m))$ be a resistance form on $V_m$ for every $m \ge 0$. A sequence of resistance forms $\{(\E_m, \l(V_m))\}_{m \ge 0}$ is called compatible if and only if $V_m \subseteq V_{m + 1}$ and $(\E_m, \l(V_m)) \le (\E_{m + 1}, \l(V_{m + 1}))$ for any $m \ge 0$.
\enddefinition

Note that if $\{(\E_m, \l(V_m))\}_{m \ge 0}$ is a compatible sequence, then, for any function $u\colon\cup_{m \ge 0} V_m \to \BbR$, the sequence $\E_m(u|_{V_m}, u|_{V_m})$ is monotonically non-decreasing. By this fact, the following definition makes sense.

\definition\label{BRF.def25}
Let $V_m$ be a finite set and let $(\E_m, \l(V_m))$ be a resistance form on $V_m$ for every $m \ge 0$. If $\S = \{(\E_m, \l(V_m))\}_{m \ge 0}$ is a compatible sequence, then we define
\[
\F_{\S} = \{u~|~u \in V_*, \lim_{m \to \infty} \E_m(u|_{V_m}, u|_{V_m}) < \infty\},
\]
where $V_* = \cup_{m \ge 0} V_m$, and for any $u, v \in \F_{\S}$,
\[
\E_{\S}(u, v) = \lim_{m \to \infty} \E_m(u|_{V_m}, v|_{V_m}).
\]
\enddefinition

\thm[Theorem~3.13 of \cite{Kig12}]\label{BRF.thm10}
Let $V_m$ be a finite set and let $(\E_m, \l(V_m))$ be a resistance form on $V_m$ for every $m \ge 0$. If $\S = \{(\E_m, \l(V_m))\}_{m \ge 0}$ is a compatible sequence, then $(\E_{\S}, \F_{\S})$ is a resistance form on $V_*$. Furthermore, let $R_{\S}$ be the associated resistance metric on $V_*$ and let $(X, R)$ be the completion of $(V_*, R_{\S})$. Then, there exists a unique resistance form $(\E, \F)$ on $X$ such that for any $u \in \F$, $u$ is continuous on $(X, R)$, $u|_{V_*} \in \F_{\S}$ and $\E(u, u) = \E_{\S}(u|_{V_*}, u|_{V_*})$. In particular, $R$ is the resistance metric associated with $(\E, \F)$.
\endthm

An important concept is the notion of trace of a resistance form. This corresponds, roughly speaking, to the restriction of a resistance form to a subset of the original domain.

\definition\label{BRF.def30}
Let $(\E, \F)$ be a resistance form on a set $X$. For any $Y \subseteq X$, define
\[
\F|_Y = \{u|_{Y}~:~u \in \F\}
\]
and
\[
\F_0(Y) = \{u~|~u \in \F,~u|_{Y} \equiv 0\}.
\]
\enddefinition

\prop[Lemma~8.2 and Theorem~8.4 of~\cite{Kig12}]\label{BRF.prop20}
Let $(\E, \F)$ be a resistance form on a set $X$ and let $Y \subseteq X$ be non-empty. Then, for any $u_* \in \F|_Y$, there exists a unique $u \in \F$ such that $u|_{Y} = u_*$ and
\[
\E(u, u) = \min\{\E(v, v)~|~ v \in \F, v|_{Y} = u_*\}.
\]
Moreover, if we denote $u = h_Y(u_*)$, then the map $h_Y: \F|_Y \to \F$ is linear. If we define $\E|_Y(u, v) = \E(h_Y(u), h_Y(v))$ for any $u, v \in \F|_Y$, then $(\E|_Y, \F|_Y)$ is a resistance form on $Y$ and the associated resistance metric $R_Y$ is the restriction onto $Y \times Y$ of the resistance metric associated with $(\E, \F)$.
\endprop

\definition\label{BRF.def40}
Let $(\E, \F)$ be a resistance form on a set $X$ and let $Y \subseteq X$ be non-empty. The map $h_Y: \F|_Y \to \F$ is called the $Y$-harmonic extension map associated with $(\E, \F)$ and $h_Y(u_*)$ is called the $Y$-harmonic function with boundary value $u_*$ associated with $(\E, \F)$. We define $\H_{(\E, \F)}(Y) = h_Y(\F|_{Y})$. The resistance form $(\E|_Y, \F|_Y)$ on $Y$ is called the trace of $(\E, \F)$ on $Y$.
\enddefinition

By~\cite[Lemma~8.5]{Kig12} and the discussion after it, the domain of a resistance form admits the orthogonal decomposition presented below.

\prop\label{BRF.prop30}
Let $(\E, \F)$ be a resistance form on a set $X$ and let $Y \subseteq X$ be non-empty. Then,
\[
\F = \H_{(\E, \F)}(Y) \oplus \F_0(Y),
\]
where $\oplus$ represents the direct sum. Moreover, for any $u \in \F$, the projection of $u$ onto $\H_{(\E, \F)}(Y)$ associated with the above direct sum is given by $h_Y(u|_Y)$ and
\[
\E(u, u) = \E(h_Y(u), h_Y(u)) + \E(u - h_Y(u), u - h_Y(u)).
\]
\endprop
Finally, Theorem~\ref{BRF.thm10} along with \cite[Theorem~3.14]{Kig12} leads to the following result.
\thm\label{BRF.thm20}
Let $(\E, \F)$ be a resistance form on a set $X$ and let $R$ be the associated resistance metric on $X$. If $\{V_m\}_{m \ge 0}$ is an increasing sequence of finite subsets of $X$, i.e. $V_m \subseteq V_{m + 1} \subseteq X$ for any $m \ge 0$, then $\S=\{(\E|_{V_m}, \l(V_m))\}_{m \ge 0}$ is a compatible sequence of resistance forms. If $A$ is the closure of $V_*$ with respect to $R$, then for any $u \in \F|_A$, $u|_{V_*} \in \F_\S$ and $\E|_A(u, u) = \E_{\S}(u|_{V_*}, u|_{V_*})$.
\endthm

\setcounter{equation}{0}

\section{Construction of resistance forms on $K$}\label{CRF}
In this section, we explain how to construct resistance forms on $K$ by means of compatible sequences in a natural way that takes into full consideration the intrinsic symmetry of $K$.
\prop\label{CRF.prop10}
$(Q_0^{\SS}, \l(V_0))$ is a resistance form on $V_0$.
\endprop
\begin{proof}
Since $V_0$ is a finite set, all properties of a resistance form (see Definition~\ref{BRF.def10}) are immediately fulfilled.
\end{proof}

\definition\label{CRF.def20}
For each $m\geq 1$, define the quadratic form $Q_m^I(\cdot, \cdot)$ on $\l(V_m)$ by
\[
Q_1^I(u, u) = \sum_{(i, j) \in B}(u(p_{ij}) - u(p_{ji}))^2
\]
for any $u \in \l(V_1)$, and by
\[
Q_m^I(u, u) = \sum_{w \in W_{m - 1}} Q_1^I(u{\circ}G_w, u{\circ}G_w)
\]
for $m \ge 2$ and any $u \in \l(V_m)$.
\enddefinition

Note that neither $Q_m^{\SS}(\cdot, \cdot)$ defined in Definition~\ref{RFF.def02} nor $Q_m^I(\cdot, \cdot)$ are resistance forms if $m \ge 1$ because $\a \in (0, 1)$. However, we show in the next lemma that any weighted combination of them actually yields a resistance form on $V_m$ for any $m\geq 1$.

\lemma\label{lemma10}
For any $m\geq 1$, let $\delta, \c_1, \ldots, \c_m$ be positive numbers. If
\[
Q(u, u) = \frac{1}{\delta}Q_m^{\SS}(u, u) + \sum_{k = 1}^m \frac1{\c_k}Q_k^I(u, u)
\]
for any $u \in \l(V_m)$, then $(Q, \l(V_m))$ is a resistance form on $V_m$.
\endlemma
\begin{proof}
Again, the conditions in Definition~\ref{BRF.def10} are fulfilled because $V_m$ is finite.
\end{proof}

As a first step to construct resistance forms on $K$, we consider compatible sequences of resistance forms on the sets $V_m$. To this purpose, we introduce the concept of \textit{matching pairs} of resistances.

\definition\label{CRF.def30}
A pair $(r, \rho) \in (0, \infty)^2$ is said to be matching if and only if
\[
\frac 53r + \rho = 1.
\]
The collection of all matching pairs of resistances will be denoted by $\MP$.
\enddefinition

\begin{figure}[H]
\centering
\begin{tikzpicture}[scale=0.4]
\coordinate[label=below:{$p_2$}] (p_1) at (0,0);
\coordinate[label=above:{$p_1$}] (p_2) at (3.5,6.0621);
\coordinate[label=below:{$p_3$}] (p_3) at (7,0);
\coordinate (p_12) at (1.5,2.59805);
\coordinate (p_21) at (2,3.46405);
\coordinate (p_13) at (3,0);
\coordinate (p_31) at (4,0);
\coordinate (p_23) at (5,3.46405);
\coordinate (p_32) at (5.5,2.59805);

\node at (3.5,-1.5) {$Q_0^\SS$};

\draw (p_1) -- (p_2) node[midway, left] {$1$} -- (p_3) node[midway, right] {$1$} -- (p_1) node[midway, above] {$1$};
\fill (p_1) circle (5pt);
\fill (p_2) circle (5pt);
\fill (p_3) circle (5pt);
\end{tikzpicture}
\hspace*{.5in} 
\begin{tikzpicture}[scale=0.4]
\coordinate (p_1) at (0,0);
\coordinate (p_2) at (3.5,6.0621);
\coordinate (p_3) at (7,0);
\coordinate (p_12) at (0.875,1.5155);
\coordinate (p_21) at (2.625,4.5466);
\coordinate (p_13) at (1.75,0);
\coordinate (p_31) at (5.25,0);
\coordinate (p_23) at (4.375,4.5466);
\coordinate (p_32) at (6.125,1.5155);

\node at (3.5,-1.5) {$E_1$};

\draw (p_1) -- (p_13) node[midway, below] {$r$} -- (p_31) node[midway, below] {$\rho$} -- (p_3) node[midway, below] {$r$};
\draw (p_1) -- (p_2) -- (p_3) (p_13) -- (p_12) (p_21) -- (p_23) (p_31) -- (p_32);
\fill (p_1) circle (5pt);
\fill (p_2) circle (5pt);
\fill (p_3) circle (5pt);
\fill (p_12) circle (5pt);
\fill (p_21) circle (5pt);
\fill (p_31) circle (5pt);
\fill (p_13) circle (5pt);
\fill (p_23) circle (5pt);
\fill (p_32) circle (5pt);
\end{tikzpicture}
\caption{Renormalization of resistances}
\label{Fig 2}
\end{figure}
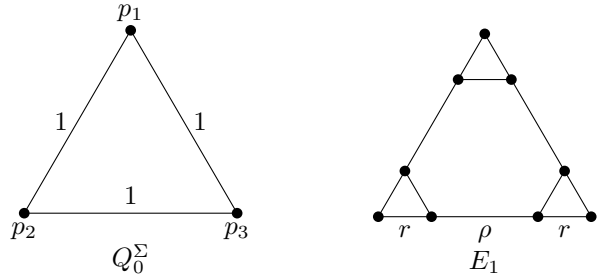

The next lemma displays the nature of the definition of matching pairs and it follows from a straightforward application of the $\Delta$-Y transform as illustrated in Figure~\ref{Fig 2}. For details on the $\Delta$-Y transform see~\cite[Lemma~2.1.15]{Kig01}.

\lemma\label{lemma20}
Let $r$ and $\rho$ be positive real numbers and define the resistance form $(E_1, \l(V_1))$ on $V_1$ by
\[
E_1(u, u) = \frac 1rQ_1^{\SS}(u, u) + \frac 1{\rho}Q_1^I(u, u)
\]
for any $u \in \l(V_1)$. Then, $(E_1, \l(V_1))$ on $V_1$ is compatible with $(Q_0^{\SS}, \l(V_0))$ on $V_0$ if and only if
$(r, \rho)$ is matching.
\endlemma

This result is the basis leading to the relationship between sequences of matching pairs and compatible sequences of resistance forms.

\thm\label{thm10}
Define 
\[
E_0(u, u) = Q^{\SS}_0(u, u)
\]
for any $u \in \l(V_0)$ and 
\[
E_m(u, u) = \frac{1}{\delta_m}Q_m^{\SS}(u, u) + \sum_{k = 1}^m \frac1{\c_{k}}Q_k^I(u, u)
\]
for any $u \in \l(V_m)$ and $m \ge 1$. Then, $\{(E_m, \l(V_m))\}_{m \ge 0}$ is a compatible sequence if and only if there exists a sequence of matching pairs $\{(r_m, \rho_m)\}_{m \ge 1}$ such that
\[
\delta_m = r_1\cdots{r_m}\quad \text{and} \quad \c_{k} = r_1\cdots{r_{k - 1}}\rho_k
\]
for any $m \ge 1$ and any $k \ge 1$.
\endthm

\begin{proof}
By definition, $\{(E_m,\ell(V_m))\}_{m\geq 0}$ is compatible if and only if $(E_m,\ell(V_m))$ is compatible with $(E_{m+1},\ell(V_{m+1}))$ for all $m\geq 0$. 

\begin{figure}[H]
\centering
\begin{tikzpicture}[scale=0.45]
\coordinate[label=below:{$p_2$}] (p_1) at (0,0);
\coordinate[label=above:{$p_1$}] (p_2) at (3.5,6.0621);
\coordinate[label=below:{$p_3$}] (p_3) at (7,0);
\coordinate (p_12) at (1.5,2.59805);
\coordinate (p_21) at (2,3.46405);
\coordinate (p_13) at (3,0);
\coordinate (p_31) at (4,0);
\coordinate (p_23) at (5,3.46405);
\coordinate (p_32) at (5.5,2.59805);

\node[color=white] at (3.5,-0.7) {$\gamma_{m+1}$};
\draw (p_1) -- (p_2) node[midway, left] {$\delta_m$} -- (p_3) node[midway, right] {$\delta_m$} -- (p_1) node[midway, above] {$\delta_m$};\fill (p_1) circle (5pt);\fill (p_2) circle (5pt);\fill (p_3) circle (5pt);
\end{tikzpicture}
\hspace*{.5in} 
\begin{tikzpicture}[scale=0.45]
\coordinate (p_1) at (0,0);
\coordinate (p_2) at (3.5,6.0621);
\coordinate (p_3) at (7,0);
\coordinate (p_12) at (0.875,1.5155);
\coordinate (p_21) at (2.625,4.5466);\coordinate (p_13) at (1.75,0);
\coordinate (p_31) at (5.25,0);
\coordinate (p_23) at (4.375,4.5466);
\coordinate (p_32) at (6.125,1.5155);

\node at (3.5,-0.675) {\small{$\gamma_{m+1}$}};
\draw (p_1) -- (p_2) -- (p_3) (p_12) -- (p_13) (p_21) -- (p_23) (p_31) -- (p_32);\draw (p_1) -- (p_13) node[midway, below] {\small{$\delta_{m+1}$}} -- (p_31) -- (p_3) node[midway, below] {\small{$\delta_{m+1}$}};\fill (p_1) circle (5pt);\fill (p_2) circle (5pt);\fill (p_3) circle (5pt);\fill (p_12) circle (5pt);\fill (p_21) circle (5pt);\fill (p_31) circle (5pt);\fill (p_13) circle (5pt);\fill (p_23) circle (5pt);\fill (p_32) circle (5pt);
\end{tikzpicture}
\caption{Renormalization of resistances}
\label{networks}
\end{figure}
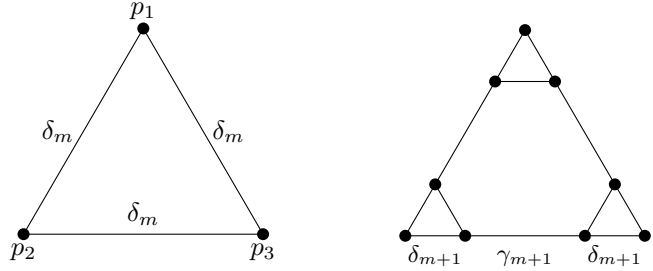

By means of the $\Delta$-Y transform, this is the case if and only if the networks in Figure~\ref{networks} are also compatible, i.e. if and only if $\frac{5}{3}\delta_{m+1}+\gamma_{m+1}=\delta_{m}$. Setting $r_m=\frac{\delta_{m+1}}{\delta_m}$, $\rho_m=\frac{\gamma_{m+1}}{\delta_m}$ and $\delta_0=r_0=\rho_0=1$, we have that $(r_m,\rho_m)$ is matching and for all $m\geq 0$,
\[
\delta_{m+1}=r_m\delta_m\qquad\text{and}\qquad\gamma_{m+1}=\rho_m\delta_m.
\]
Applying these equalities recursively leads to the desired statement.
\end{proof}

\notation
We denote by $\MP^{\BbN}$ the collection of sequences of matching pairs of resistances, i.e.
\[
\MP^{\BbN} = \big\{\{(r_m, \rho_m)\}_{m \ge 1}~|~r_m, \rho_m \in (0, \infty),~\frac 53r_m + \rho_m = 1\text{ for any }m \ge 1\big\}.
\]
\endnotation

\definition\label{RNS.def20}
 Let $\R = \{(r_m, \rho_m)\}_{m \ge 1} \in \MP^{\BbN}$ and define the quadratic form $E_{\R, m}$ on $\l(V_m)$ to be $E_m$ as given in Theorem~\ref{thm10}. Moreover, define 
\[
\widehat{\F}_{\R} = \{u~|~u \in \l(V_*), \lim_{m \to \infty} E_{\R, m}(u|_{V_m}, u|_{V_m}) < \infty\}
\]
 and 
\[
\widehat{\E}_{\R}(u, v) = \lim_{m \to \infty} E_{\R, m}(u|_{V_m}, v|_{V_m})
\]
 for any $u, v \in \widehat{\F}_{\R}$.
\enddefinition

In view of Theorem~\ref{thm10} and Theorem~\ref{BRF.thm10}, $(\widehat{\E}_{\R}, \widehat{\F}_{\R})$ is a resistance form on $\l(V_*)$ for any $\R = \{(r_m, \rho_m)\}_{m \ge 1} \in \MP^{\BbN}$. Note that if $\{(r_m, \rho_m)\}_{m \ge 1} \in \MP^{\BbN}$, then $r_m \le \frac 35$ for any $m \ge 1$ and hence $\delta_m \le (\frac 35)^m$ and $\c_m \le (\frac 35)^{m - 1}$.

\lemma\label{RNS.lemma00}
Let $\R = \{(r_m, \rho_m)\}_{m \ge 1} \in \MP^{\BbN}$ and let $R_m$ denote the resistance metric on $V_m$ associated with $(E_{\R, m}, \l(V_m))$. Then, $\diam{V_m, R_m} \le 4$ for any $m \ge 1$, where $\diam{X, d}$ is the diameter of the metric space $(X, d)$ given by $\sup_{x, y \in X} d(x, y)$. In particular
\[
|u(x) - u(y)|^2 \le 4E_{\R,m}(u, u)
\]
for any $x,y\in V_m$ and $u\in\ell(V_m)$.
\endlemma
\demo
Let $q = G_{w_1\ldots{w_m}}(p_i)$. Define $q_k = G_{w_1\ldots{w_k}}(p_{w_k})$ for $k = 1, \ldots, m$ and set $q_{m + 1} = q$. Since $G_{w_k}(p_{w_k}) = p_{w_k}$, we have that $q_k = G_{w_1\ldots{w_{k - 1}}}(p_{w_k})$ and in particular $q_1 = G_{w_1}(p_{w_1}) = p_{w_1}$. Since $\{(E_m, V_m)\}_{m \ge 0}$ is compatible, it holds that
\begin{multline*}
R_m(q_k, q_{k + 1}) = R_m(G_{w_1\ldots{w_k}}(p_{w_k}), G_{w_1\ldots{w_k}}(p_{w_{k + 1}}))\\ 
= R_k(G_{w_1\ldots{w_k}}(p_{w_k}), G_{w_1\ldots{w_k}}(p_{w_{k + 1}})) \le \delta_k.
\end{multline*}
Therefore,
\[
R_m(p_{w_1}, q) = R_m(q_1, q_{m + 1}) \le \sum_{k = 1}^m R_m(q_k, q_{k + 1}) \le \sum_{k = 1}^m \delta_k \le \sum_{k = 1}^{\infty} \bigg(\frac 35\bigg)^k =\frac 32.
\]
Thus if $x = G_{w_1\ldots{w_m}}(p_i)$ and $y = G_{v_1\ldots{v_m}}(p_j)$, then
\[
R_m(x, y) = R_m(x, p_{w_1}) + R_m(p_{w_1}, p_{v_1}) + R_m(p_{v_1}, y) \le \frac 32 + 1 +  \frac 32 = 4.
\]
\enddemo

 The resistance form $(\widehat{\E}_\R,\widehat{\F}_\R)$ possesses every symmetry and similarity (which is inhomogeneous with respect to $m$) required by the definition of completely symmetric resistance forms, although a function $u$ in the domain $\widehat{\F}_\R$ is not a function on $K$ but on $V_*$.

\lemma\label{RNS.lemma03}
Let $\R = \{(r_m, \rho_m)\}_{m \ge 1} \in \MP^{\BbN}$. If $\R^{(n)} = \{(r_{n + m}, \rho_{n + m})\}_{m \ge 1}$, then $u{\circ}G_w \in \hF_{\R^{(m)}}$ for any $m \ge 1$, $u \in \hF_{\R}$ and $w \in W_m$. Moreover, 
\[
 \hE_{\R}(u, v) = \sum_{w \in W_m} \frac 1{\delta_m}\hE_{\R^{(m)}}(u{\circ}G_w, v{\circ}G_w) + \sum_{k = 1}^m \frac 1{\c_{k}}Q_k^I(u, v)
\]
for any $u, v \in \hF_{\R}$ and $m \ge 1$.
\endlemma

\begin{proof}
Let $n,m\geq 1$. For any $u\in\widehat{\F}_\R$ it holds that
\begin{align*}
\widehat{\E}_{\R,n+m}(u,u)&=\frac{1}{\delta_{n+m}}Q_{n+m}^\SS(u,u)+\sum_{k=m+1}^{n+m}\frac{1}{\gamma_k}Q_k^I(u,u)+\sum_{k=1}^{m}\frac{1}{\gamma_k}Q_k^I(u,u)\\
&=\sum_{w\in W_m}\frac{1}{\delta_m}\widehat{\E}_{\R^{(m)},n}(u{\circ}G_w,u{\circ}G_w)+\sum_{k=1}^m\frac{1}{\gamma_k}Q_k^I(u,u).
\end{align*}
Letting $n\to\infty$ in both sides of the equality leads to the desired result, which implies that $u{\circ}G_w\in\widehat{\F}_{\R^{(m)}}$ for any $u\in\widehat{\F}_\R$ and $w\in W_m$.
\end{proof}

\lemma\label{RNS.lemma100}
Let $\R = \{(r_m, \rho_m)\}_{m \ge 1} \in \MP^{\BbN}$ and let $R$ be the resistance metric on $V_*$ associated with $(\hE_{\R}, \hF_{\R})$. If $\overline{V_*}$ is the completion of $V_*$ with respect to $R$, then the identity map $\iota: V_* \to V_*$ is extended to a homeomorphism from $(\overline{V_*}, R)$ to $(\SS, d_E)$.
\endlemma

By the above lemma and Theorem~\ref{BRF.thm10}, the resistance form $(\hE_{\R}, \hF_{\R})$ is naturally regarded as a resistance form on $\SS$ and $\hF_{\R}$ is thought of as a subset of $C(\SS)$.\par

\demo
Let $w \in W_*$ and let $x, y \in V_*$. Set $p = G_w(x)$ and $q = G_w(y)$. By Lemma~\ref{RNS.lemma00} and Lemma~\ref{RNS.lemma03}, 
\[
\frac{|u(p) - u(q)|^2}{\E_{\R}(u, u)} \le \delta_m\frac{|u(G_w(x)) - u(G_w(y))|^2}{\E_{\R^{(m)}}(u{\circ}G_w, u{\circ}G_w)} \le 4\delta_m
\]
holds for any $u \in \hF_{\R}$. Thus, $\diam{G_w(V_*), R} \le 4\delta_m$. Let $\{x_n\}_{n \ge 1}$ be a Cauchy sequence in $V_*$ with respect to $d_E$. Then, there exists $x \in \SS$ such that $d_E(x, x_n) \to 0$ as $n \to \infty$. If $x = \iota(\omega_1\omega_2\ldots)$, then $x_n \in G_{\word{\omega}m}(V_*)$ for sufficiently large $n$ and therefore $R(x_k, x_l) \le 4\delta_m$ for sufficiently large $k$ and $l$. Hence, $\{x_n\}_{n \ge 1}$ is a Cauchy sequence in $(V_*, R)$ as well. \par
On the other hand, if $w, v \in W_m$ and $w \neq v$, by (RF3), there exists $u \in \hF_{\R}$ such that $u|_{G_w(V_*)} \equiv 1$ and $u|_{G_v(V_*)} \equiv 0$. For any $x \in G_w(V_*)$ and $y \in G_v(V_*)$, we thus have that
\[
1 = |u(x) - u(y)|^2 \le \hE_{\R}(u, u)R(x, y),
\]
which shows that $\inf\{R(x, y)~|~x \in G_w(V_*), y \in G_v(V_*)\} > 0$. Hence, 
\[
\min\Big\{\inf\{R(x, y)~|~x \in G_w(V_*), y \in G_v(V_*)\}~\big|~w, v \in W_m, w \neq v\Big\} > 0.
\]
If $\{y_n\}_{n \ge 1}$ is a Cauchy sequence in $(V_*, R)$, then for any $m \ge 0$, there exists $w \in W_m$ such that $y_n \in G_w(V_*)$ for sufficiently large $n$. Thus, $\{y_n\}_{n \ge 1}$ is a Cauchy sequence in $(V_*, d_E)$ as well.\\
Consequently, the identity map from $V_*$ to $V_*$ is extended to a homeomorphism between $\overline{V_*}$ and $\SS$.
\enddemo

In order to obtain a resistance form on $K$ from a sequence of matching pairs of resistances, we need to replace $Q^I_k(u, u)$ by a sum of $H^1$-inner products on the line segments $e_{ij}^w$. The bilinear form that arises in this way preserves many properties of the former and it is defined for functions in $\widetilde{\F}$.

\definition\label{RNS.def500}
Let $\R = \{(r_m, \rho_m)\}_{m \ge 1} \in \MP^{\BbN}$. For each $m\geq 1$, define
\[
\E_{\R, m}(u, v) = \frac 1{\delta_m}Q_m^{\SS}(u, v) + \sum_{k = 1}^m \frac 1{\gamma_k}\D_k^I(u, v)
\]
for any $u, v \in \widetilde{\F}$, where $\delta_m = r_1\cdots{r_m}$ and $\gamma_m=\delta_{m- 1}\rho_m$.
\enddefinition

We start by establishing some relations between the forms $E_{\R,m}$ and $\E_{\R,m}$.

\lemma\label{RNS.lemma05}
For any $u \in \tF$ and $m\geq 1$,
\begin{equation}\label{RNS.eq100}
Q_m^I(u, u) \le \D^I_m(u, u)
\end{equation}
and
\begin{equation}\label{RNS.eq105}
E_{\R, m}(u, u) \le \E_{\R, m}(u, u).
\end{equation}
\endlemma
\demo
For any $u \in H^1([0, 1])$,
\[
(u(1) - u(0))^2 \le \int_0^1 \bigg(\frac{du}{dx}\bigg)^2dx.
\]
Applying this to every $e_{ij}^w$, we obtain~\eqref{RNS.eq100}. Consequently,
\begin{multline*}
E_{\R, m}(u, u) = \frac{1}{\delta_m}Q_m^{\SS}(u, u) + \sum_{k = 1}^m \frac 1{\c_k}Q_k^I(u, u)
\\ \le \frac{1}{\delta_m}Q_m^{\SS}(u, u) + \sum_{k = 1}^m \frac 1{\c_k}\D_k^I(u, u) = \E_{\R, m}(u, u).
\end{multline*}
\enddemo

\lemma\label{RNS.lemma10}
Assume that $(r, \rho)\in \MP$. Then, for any $u \in \tF$,
\begin{equation}\label{RNS.eq110}
Q_0^{\SS}(u, u) \le \frac 1{r}Q^{\SS}_1(u, u) + \frac 1{\rho}\D^I_1(u, u).
\end{equation}
\endlemma

\demo
By Lemma~\ref{lemma20}, we see that
\[
Q_0^{\SS}(u, u) \le \frac 1{r}Q^{\SS}_1(u, u) + \frac 1{\rho}Q^I_1(u, u).
\]
Combining this with Lemma~\ref{RNS.lemma05}, we obtain \eqref{RNS.eq110}.
\enddemo

We can now use these properties in order to show that for any $u\in\widetilde{\F}$, the sequence $\{\E_{\R,m}(u,u)\}_{m\geq 1}$ is monotonically non-decreasing.

\lemma\label{RNS.lemma20}
If $\{(r_m, \rho_m)\}_{m \ge 1} \in \MP^{\BbN}$, then
\[
\E_{\R, m}(u, u) \le \E_{\R, m + 1}(u, u)
\]
for any $u \in \tF$.
\endlemma
\demo
By Lemma~\ref{RNS.lemma10}, it follows that
\begin{multline*}
\sum_{w \in W_m}Q_0^{\SS}(u{\circ}G_w, u{\circ}G_w) \\
\le \sum_{w \in W_m}\Big(\frac 1{r_{m + 1}} Q_1^{\SS}(u{\circ}G_w, u{\circ}G_w) + \frac{1}{\rho_{m + 1}}\D_1^{I}(u{\circ}G_w, u{\circ}G_w)\Big).
\end{multline*}
Multiplying by $(\delta_m)^{-1}$ and adding $\sum_{k = 1}^m \frac 1{\delta_{k - 1}\rho_k}\D_k^I(u, u)$ on both sides of the inequality, we verify the desired statement.
\enddemo

In view of this lemma, $\{\E_{\R, m}(u, u)\}_{m\geq 1}$ converges to a non-negative real number or infinity as $m \to \infty$. Therefore, the following definition makes sense.

\definition\label{RNS.def10}
Let $\R = \{(r_m, \rho_m)\}_{m \ge 1} \in \MP^{\BbN}$. Define
\[
\F_{\R} = \{u~|~u \in \widetilde{\F}, \lim_{m \to \infty} \E_{\R, m}(u, u) < \infty\}
\]
and
\[
\E_{\R}(u, v) = \lim_{m \to \infty} \E_{\R, m}(u, v)
\]
for any $u, v \in \F_{\R}$. 
\enddefinition

 The next theorem is the main result of this section. It shows that resistance forms on $K$ constructed from a sequence of matching pairs $\R$ are completely symmetric resistance forms. In addition, it provides an explicit expression of their corresponding resolution.

\thm\label{thm20}
Let $\R = \{(r_m, \rho_m)\}_{m \ge 1} \in \MP^{\BbN}$. Then, $(\E_{\R}, \F_{\R}) \in \RF_S$. More precisely, if $\R^{(n)} = \{(r_{n + m}, \rho_{n + m})\}_{m\geq 1} \in \MP^{\BbN}$ for any $n \ge 0$, then the resolution of $(\E_{\R}, \F_{\R})$ is given by $\{((\delta_m)^{-1}\E_{\R^{(m)}}, \F_{\R^{(m)}}, \c_m)\}_{m \ge 0}$, where $\delta_m = r_1{\cdots}r_m$ and $\c_m = \delta_{m - 1}\rho_m$ for any $m \ge 0$.
\endthm

In order to show this theorem, we need several lemmas.

\lemma\label{RNS.lemma30}
Let $\R = \{(r_m, \rho_m)\}_{m \ge 1} \in \MP^{\BbN}$. For any $x \neq y \in K$, there exists $u \in \tF_\infty$ such that $u(x) \neq u(y)$.
\endlemma

\demo
If either $x$ or $y$ belong to $\sd{K}{\SS}$, for instance $x \in \sd{K}{\SS}$, then there exists $w \in W_*$ such that $x \in G_w(\sd{e_{ij}}{\{p_{ij}, p_{ji}\}})$. In this case, there exists $u|_{e_{ij}^w} \in H^1(e_{ij}^w)$ such that $u|_{e_{ij}^w}(x) = 1$ and $u|_{e_{ij}^w}(G_w(p_{ij})) = u|_{e_{ij}^w}(G_w(p_{ji})) = 0$. Letting $u(z) = 0$ for any $z \in \sd{K}{e_{ij}^w}$, we obtain the desired function $u \in \tF_{|w| + 1}$.\par
If $x, y \in \SS$, then there exist $m \ge 1$, $w \in W_m$ and $v \in W_m$ such that $x \in G_w(K)$, $y \in G_v(K)$ and $G_w(K) \cap G_v(K) = \emptyset$. Now, there is a function $u \in \tF_m$ such that $u|_{G_w(K)} = 1$ and $u|_{G_v(K)} = 0$.
\enddemo

The following lemma is straightforward from the definition of $\E_{\R, m}$, $\F_{\R}$ and $\E_{\R}$.

\lemma\label{RNS.lemma40}
Let $\R = \{(r_m, \rho_m)\}_{m \ge 1} \in \MP^{\BbN}$ and let $\R^{(n)} = \{(r_{m + n}, \rho_{m + n})\}_{m \ge 1}$ for any $n \ge 0$. Then,
\begin{multline*}
\F_{\R} = \{u~|~u\colon K \to \BbR,~u{\circ}G_w \in \F_{\R^{(n)}}\text{ for any }w \in W_n,\\
 u|_{e_{ij}^w} \in H^1(e_{ij}^w)\text{ for any }(w, (i, j)) \in \big(\cup_{k = 0}^{n - 1}W_k\big) \times B\}
\end{multline*}
and for any $u \in \F_{\R}$,
\[
\E_{\R}(u, u) = \frac 1{\delta_n}\sum_{w \in W_n} \E_{\R^{(n)}}(u{\circ}G_w, u{\circ}G_w) + \sum_{k = 1}^n \frac 1{\c_k}\D_k^I(u, u).
\]
\endlemma

\lemma\label{RNS.lemma50}
Let $\R = \{(r_m, \rho_m)\}_{m \ge 1} \in \MP^{\BbN}$. If $w \in W_m$, then
\begin{equation}\label{RNS.eq200}
|u(x) - u(y)|^2 \le 16\delta_{m}\E_{\R}(u, u)
\end{equation}
for any $u\in\F_\R$ and $x, y \in G_w(K)$.
\endlemma

\demo
First we show the case when $m = 0$ and $w = \emptyset$, namely 
\begin{equation}\label{RNS.eq210}
|u(x) - u(y)|^2 \le 16\E_{\R}(u, u)
\end{equation}
 for any $u \in \F_{\R}$ and $x, y \in K$. Let $x, y \in B_*:=\bigcup_{(w,(i,j))\in W_*\times B}e_{ij}^w$. Then, $x \in G_{w}(e_{ij}), y \in G_v(e_{kl})$ for some $w, v \in W_*$ and $(i, j),(k,l)\in B$. Set $p = G_w(p_i)$ and $q = G_v(p_k)$. Since $\c_n=\delta_{n-1}\rho_n\le 1$ for any $n\geq 1$, it follows that
\[
|u(x) - u(p)|^2 \le \c_{|w| + 1}\frac{1}{\c_{|w| + 1}}\D_{e_{ij}^w}(u, u) \le \c_{|w| + 1}\E_{\R, |w| + 1}(u, u)  \le \E_{\R}(u, u).
\]
In the same way, we obtain $|u(y) - u(q)|^2 \le \E_{\R}(u, u)$. Setting $m = \max\{|w|, |v|\}$, Lemma~\ref{RNS.lemma00} and \eqref{RNS.eq105} yield
\[
|u(p) - u(q)|^2 \le 4E_{\R, m}(u, u) \le 4\E_{\R, m}(u, u) \le 4\E_{\R}(u, u).
\]
Combining these inequalities, we have
\[
|u(x) - u(y)|^2 \le (|u(x) - u(p)| + |u(p) - u(q)| + |u(q) - u(y)|)^2 \le 16\E_{\R}(u, u).
\]
Since $\tF \subseteq C(K)$ and $B_*$ is dense in $K$ with respect to the Euclidean metric, \eqref{RNS.eq210} holds for any $x, y \in K$.

Consider now $w \in W_m$ with $m\geq 1$, and set $x = G_w(x')$ and $y = G_w(y')$. For any $u \in \F_{\R}$, Lemma~\ref{RNS.lemma40} implies that $u{\circ}G_w \in \F_{\R^{(m)}}$. Applying \eqref{RNS.eq210} to $(\E_{\R^{(m)}}, \F_{\R^{(m)}})$ and using again Lemma~\ref{RNS.lemma40}, we see that
\begin{multline*}
|u(x) -u(y)|^2 = |u(G_w(x')) - u(G_w(y'))|^2
 \le 16\E_{\R^{(m)}}(u\circ{G_w}, u\circ{G_w}) \le 16\delta_m\E_{\R}(u, u).
\end{multline*}
\enddemo

\demo[Proof of Theorem~\ref{thm20}]
We start by showing that $(\E_{\R}, \F_{\R})$ is a resistance form on $K$.\\
(RF1): By definition, $\F_{\R} \subseteq C(K)$ and $\E_{\R}$ is a non-negative quadratic form on $\F_{\R}$. Moreover, if $\E_{\R}(u, u) = 0$, then $\E_{\R, m}(u, u) = 0$ for any $m \ge 0$. This implies that $u$ is constant on $e_{ij}^w$ and $G_w(V_0)$ for any $(w,(i, j))\in W_m\times B$. Therefore, $u$ is constant on $K$ and (RF1) holds.\\
(RF2): It suffices to prove that $(\F_{\R, 0}, \E_{\R})$ is complete, where $\F_{\R, 0} = \{u | u \in \F_{\R}, u(p_1) = 0\}$. Let $\{u_n\}_{n \ge 1}$ be a Cauchy sequence in $(\F_{\R, 0}, \E_{\R})$. By~\eqref{RNS.eq200},
\[
|u_n(x) - u_m(x)|^2 = |(u_n - u_m)(p_1) - (u_n - u_m)(x)|^2 \le 16\E_{\R}(u_n - u_m, u_n - u_m).
\]
This implies that $\{u_n\}_{n\geq 1}$ converges uniformly on $K$ as $n \to \infty$. Let $u$ be its limit. Then $u_n|_{e_{ij}^w}$ converges to $u|_{e_{ij}^w}$ in the sense of $H^1(e_{ij}^w)$ and hence $u|_{e_{ij}^w} \in H^1(e_{ij}^w)$. If $m \ge n$, 
\[
\E_{\R, k}(u_n - u_m, u_n - u_m) \le \E_{\R}(u_n - u_m, u_n - u_m) \le \sup_{m \ge n} \E_{\R}(u_n - u_m, u_n - u_m).
\]
Letting first $m \to \infty$ and afterwards $k \to \infty$, we see that $u \in \F_{\R}$ and
\[
\E_{\R}(u_n -u, u_n - u) \le \sup_{m \ge n} \E_{\R}(u_n - u_m, u_n - u_m).
\]
Letting $n \to \infty$, we finally verify $\E_{\R}(u_n - u, u_n - u) \to 0$ as $n \to \infty$. Thus, $(\F_{\R, 0}, \E_{\R})$ is complete.\\
(RF3) follows from Lemma~\ref{RNS.lemma30}.\\
(RF4) is immediate by Lemma~\ref{RNS.lemma50}.\\
(RF5): Note that for any $u \in \F_{\R}$ and any $m\geq 1$,
\begin{equation}\label{RNS.eq300}
Q_m^{\SS}(\ou, \ou) \le Q_m^{\SS}(u, u)\quad\text{and}\quad \D_m^I(\ou, \ou) \le \D_m^I(u, u).
\end{equation}
This implies that $\E_{\R, m}(\ou, \ou) \le \E_{\R, m}(u, u)$ for any $m \ge 1$, hence $\ou \in \F_{\R}$ and $\E_{\R}(\ou, \ou) \le \E_{\R}(u, u)$.

Thus we have shown that $(\E_{\R}, \F_{\R})$ is a resistance form on $K$. Let us prove next that the identity map from $(K, d_E)$ to $(K, R)$ is continuous. Assume that $d_E(x_n, x) \to 0$ as $n \to \infty$. For the moment, we consider the following two cases (I) and (II):\\
(I) There exists $(w,(i, j))\in W_*\times B$ such that $\{x_n\}_{n \ge 1} \subseteq e_{ij}^w$ and $x \in e_{ij}^w$.\\
(II) There exists $\{w(n)\}_{n \ge 1} \subseteq W_*$ such that $x, x_n \in G_{w(n)}(K)$ and $\lim_{n \to \infty} |w(n)| = \infty$.\\
Assume (I). Since for any $u\in\F_\R$
\[
\frac{|u(x_n) - u(x)|^2}{\E_{\R}(u, u)} \le \frac{|u(x_n) - u(x)|^2}{{\c_{|w|}}^{-1}\D_{e_{ij}^w}(u, u)} \le \c_{|w|}\frac{d_E(x_n, x)}{d_E(G_w(p_{ij}), G_w(p_{ji}))},
\]
it follows that $R(x_n, x) \le \c_md_E(x_n, x)/d_E(G_w(p_{ij}), G_w(p_{ji}))$ and hence $R(x_n, x) \to 0$ as $n \to \infty$.\\
Assume (II). Then Lemma~\ref{RNS.lemma50} yields $R(x_n, x) \le 16\delta_{|w(n)|}$, which immediately implies that $R(x_n, x) \to 0$ as $n \to \infty$.\\
Let us now consider general cases. If $x \in \sd{K}{B_*}$, then there exists $w_1w_2\ldots \in S^{\BbN}$ such that $x = \bigcap_{m \ge 1} G_{\word wm}(K)$ and $x$ belongs to the interior of $G_{\word wm}(K)$ for any $m \ge 1$. Thus, if $d_E(x_n, x) \to 0$ as $n \to 0$, then we have case (II). If $x$ belongs to the interior of $e_{ij}^w$ for some $(w,(i,j)) \in W_*\times B$, then we have case (I). Finally, if $x = G_w(p_{ij})$ and $d_E(x_n, x) \to 0$ as $n \to 0$, then we can decompose $\{x_n\}_{n \ge 1}$ into $\{x_n~|~x_n \in e_{ij}^w\}$ and $\{x_n~|~x_n \in \sd{K}{e_{ij}^w}\}$ (either one may be empty). Applying case (I) and case (II) to the first part and the second part respectively, we verify that $R(x_n, x) \to 0$ as $n \to \infty$. Thus, we have shown that the identity from $(K, d_E)$ to $(K, R)$ is continuous. Since $(K, d_E)$ is compact, so is $(K, R)$ and the inverse is continuous as well. Therefore, $R$ gives the same topology as $d_E$. Notice that by definition, $\E_{\R,m}(u,u)=\E_{\R,m}(u{\circ}\vp,u{\circ}\vp)$ for any $\vp\in\G_K$ and hence the same holds for $\E_\R$.\par
Finally, applying Lemma~\ref{RNS.lemma40}, we conclude that $(\E_{\R}, \F_{\R}) \in \RF_S$ and its resolution is $\{((\delta_m)^{-1}\E_{\R^{(m)}}, \F_{\R^{(m)}}, \c_m)\}_{m \ge 0}$.
\enddemo

\setcounter{equation}{0}

\section{Identification of $\RF_S$ with the resistance forms \\
from matching pairs}\label{IRF}

In the previous section, we proved that any resistance form $(\E_\R,\F_\R)$ derived from a sequence of matching pairs $\R$ is completely symmetric. This section focuses on the converse statement by proving in Theorem~\ref{IRF.thm20} that, up to multiplication by a constant, any completely symmetric resistance form can be obtained from a sequence of matching pairs.\par
First of all, notice that since $(\E, \F) \in \RF_S^{(0)}$ is symmetric and $V_0$ has only three points, one immediately arrives to the following fact.

\lemma\label{IRF.lemma10}
 For any $(\E, \F) \in \RF_S^{(0)}$, there exists $r_0 > 0$ such that 
\[
\E_{V_0}(u, v) = \frac 1{r_0}Q_0^{\SS}(u, v)
\]
for any $u, v \in \l(V_0)$, where $\E_{V_0}$ is the trace of $(\E, \F)$ on $V_0$.
\endlemma

See Proposition~\ref{BRF.prop20} for the definition of trace of a resistance form.

\demo
Since $(\E, \F) \in \RF_S^{(0)}$, the trace $\E_{V_0}$ has the same symmetry as the equilateral triangle $p_1p_2p_3$. Therefore, $\E_{V_0}$ must be a constant multiple of $Q_0^{\SS}$.
\enddemo

Resistance forms whose trace on $V_0$ coincides with $Q_0^\SS$ will play a special role in the forthcoming discussion.

\definition\label{IRF.def10}
For $(\E, \F) \in \RF_S^{(0)}$, define $r_0(\E, \F)$ to be the constant $r_0$ given in Lemma~\ref{IRF.lemma10}. Furthermore, define
\[
\RF_S^N = \{(\E, \F)~|~(\E, \F) \in \RF_S,~r_0(\E, \F) = 1\}.
\]
\enddefinition
The superscript ``N'' in $\RF_S^N$ represents the word ``normalized'' since $r_0(\E, \F) = 1$.
\lemma\label{IRF.lemma20}
Let $(\E, \F) \in \RF_S$ and let $\{(\E_m, \F_m, \eta_m)\}_{m \ge 0} \subseteq \RF_S^{(0)} \times (0, \infty)$ be the resolution of $(\E, \F)$. If $\delta_m = r_0(\E_m, \F_m)$ for $m \ge 0$, then $(\delta_m/\delta_{m - 1}, \eta_m/\delta_{m - 1}) \in \MP$.
\endlemma

\demo
From Definition~\ref{RFF.def10}, we know that $\F_m=\{u{\circ}G_i~|~u \in \F_{m - 1}\}$ for any $i \in S$ and $\{u|_{e_{ij}}: u \in \F_{m  - 1}\} = H^1(e_{ij})$ for any $(i, j) \in B$. Therefore, equality~\eqref{RFF.eq10} yields
\[
\E_{m - 1}|_{V_1}(u, v) = \frac 1{\delta_{m}}Q_1^{\SS}(u, v) + \frac 1{\eta_m}Q^I_1(u, v)
\]
for any $u, v \in \F_{m - 1}|_{V_1} = \l(V_1)$, where $(\E_{m - 1}|_{V_1}, \F_{m - 1}|_{V_1})$ is the trace of $(\E_{m - 1}, \F_{m - 1})$ on $V_1$ (see Definition~\ref{BRF.def30}).  On the other hand, 
\[
\E_{m - 1}|_{V_0}(u, v) = \frac 1{\delta_{m - 1}}Q_0^{\SS}(u, v)
\]
for any $u, v \in \F_{m - 1}|_{V_0} = \l(V_0)$. Since $(\E_{m - 1}|_{V_1}, \l(V_1))$ and $(\E_{m - 1}|_{V_0}, \l(V_0))$ are compatible, Lemma~\ref{lemma20} shows that $(\delta_m/\delta_{m - 1}, \eta_m/\delta_{m - 1}) \in \MP$.
\enddemo

Due to the above lemma, it is possible to associate a sequence of matching pairs to any $(\E, \F) \in \RF_S$.

\definition\label{IRF.def20}
Let $(\E, \F) \in \RF_S$ and let $\{(\E_m, \F_m, \eta_m)\}_{m \ge 0} \subseteq \RF_S^{(0)} \times (0, \infty)$ be the resolution of $(\E, \F)$. We define $\R_{(\E, \F)} \in \MP^{\BbN}$ by $\R_{(\E, \F)} =\{ (\delta_m/\delta_{m - 1}, \eta_m/\delta_{m - 1})\}_{m \ge 1}$, where $\delta_m = r_0(\E_m, \F_m)$.
\enddefinition
We show next that for each $(\E,\F)\in\RF_S$ and $m\geq 1$, multiplied by $r_0(\E, \F)$, its trace $\E|_{V_m}$ on $V_m$ coincides with the resistance form introduced in Definition~\ref{RNS.def20} associated with the sequence of matching pairs $\R_{(\E,\F)}$.
\lemma\label{IRF.lemma30}
For any $(\E, \F) \in \RF_S$ and any $m \ge 1$, $r_0(\E, \F)\E|_{V_m} = E_{\R_{(\E, \F)}, m}$.
\endlemma

\demo
Let $\{(\E_m, \F_m, \eta_m)\}_{m \ge 0}$ be the resolution of $(\E, \F)$ and set $\delta_m = r_0(\E_m, \F_m)$. By Proposition~\ref{RFF.prop10}, $\F_m=\{u{\circ} G_w~|~u \in \F\}$, and $\F|_{e_{ij}^w} = H^1(e_{ij}^w)$ for any $m \ge 1$ and any $(w,(i, j))\in W_m\times B$. Hence, it follows from~\eqref{RFF.eq20} that
\begin{equation}\label{IRF.eq40}
\E|_{V_m}(u, v) = \sum_{w \in W_m} \frac 1{\delta_m}Q_0^{\SS}(u{\circ}G_w, u{\circ}G_w) + \sum_{k = 1}^m \frac 1{\eta_k} Q_k^I(u, v) = \frac 1{\delta_0}E_{\R, m}(u, v)
\end{equation}
for any $u, v \in \l(V_m)$.
\enddemo

\lemma\label{IRF.lemma40}
Let $(\E, \F) \in \RF_S$ and let $u_* \in \F|_{\SS}$. If $u$ is the $\SS$-harmonic function with respect to $(\E, \F)$ with boundary value $u_*$, then 
\[
u|_{e_{ij}^w}((1 - t)G_w(p_{ij}) + tG_w(p_{ji})) = (1 - t)u_*(G_w(p_{ij})) + tu_*(G_w(p_{ji}))
\]
for any $t \in [0, 1]$ and any $(w, (i, j)) \in W_* \times B$.
\endlemma

\demo
Since the restriction of a $\SS$-harmonic function to a line segment $e^w_{ij}$ is a harmonic function with respect to the Dirichlet integral, it must be an affine function.
\enddemo

The harmonic functions determined in the previous lemma provide the definition of the trace of $(\E,\F)$ on $\SS$. We will denote the subspace of $\SS$-harmonic functions by $\H_{(\E,\F)}(\SS)$ and its counterpart by $\F(\SS)$. The latter domain is characterized in the following lemma.

\lemma\label{IRF.lemma50}
Let $(\E, \F) \in \RF_S$ and let $\{(\E_m, \F_m, \eta_m)\}_{m \ge 0} \subseteq \RF_S^{(0)} \times (0, \infty)$ be the resolution of $(\E, \F)$.  If $\F(\SS) = \{u~|~u \in \F, u|_{\SS} = 0\}$, then
\begin{multline}\label{IRF.eq50}
\F(\SS) = \{u~|~u\colon K \to \BbR,~u|_{e_{ij}^w} \in H^1(e_{ij}^w)\,\text{ for any }(w, (i, j)) \in W_* \times B, \\
u|_{\SS} \equiv 0,~\sum_{k = 1}^{\infty} \frac 1{\eta_k}\D_k^I(u,u) < +\infty\}.
\end{multline}
Moreover, for any $u \in \F(\SS)$,
\begin{equation}\label{IRF.eq60}
\E(u, u) = \sum_{k = 1}^{\infty} \frac 1{\eta_k}\D_k^I(u, u).
\end{equation}
\endlemma

\demo
By Proposition~\ref{RFF.prop10}, if $u \in \F(\SS)$, then $u$ belongs to the set on the right-hand side of~\eqref{IRF.eq50}. Conversely, suppose that $u$ belongs to the set on the right-hand side of~\eqref{IRF.eq50}. Define $u_n : K \to \BbR$ as
\[
u_n(x) = \begin{cases}
u(x)\quad&\text{if }x \in \bigcup\limits_{(w,(i,j))\in \big(\cup_{k = 0}^{n - 1}W_k\big)\times B} e_{ij}^w,\\
0\quad&\text{otherwise.}
\end{cases}
\]
In view of~\eqref{RFF.eq15}, $u_n \in \F$, and if $m \ge n$, it follows from~\eqref{RFF.eq20} that
\[
\E(u_n - u_m, u_n - u_m) = \sum_{k = n + 1}^m \frac 1{\eta_k}\D_k^I(u,u).
\]
Because $\sum_{k = 1}^{\infty} \frac 1{\eta_k}\D_k^I(u, u) < \infty$, the sequence $\{u_n\}_{n \ge 1}$ is Cauchy in $(\E, \F_{p_1})$, where $\F_{p_1} = \{u~|~u \in \F, u(p_1) = 0\}$ and since $(\E, \F_{p_1})$ is complete, there exists $\tilde{u} \in \F_{p_1}$ such that $\E(\tilde{u} - u_n, \tilde{u} - u_n) \to 0$ as $n \to \infty$. Moreover, for any $x \in K$,
\[
|u_n(x) - \tilde{u}(x)|^2 = |u_n(x) - \tilde{u}(x) - (u_n(p_1) - \tilde{u}(p_1))|^2 \le \E(u_n - \tilde{u}, u_n - \tilde{u})R(x, p_1),
\]
where $R(\cdot, \cdot)$ is the resistance metric associated with $(\E, \F)$. This implies that $u_n(x) \to \tilde{u}(x)$ as $n \to \infty$. On the other hand, for any $x \in K$, $u_n(x) \to u(x)$ as $n \to \infty$. Hence, $u = \tilde{u} \in \F$ and
\[
\E(u, u) = \lim_{n \to \infty}\E(u_n, u_n) = \sum_{k = 1}^{\infty} \frac 1{\eta_k}\D_k^I(u, u).
\]
\enddemo

We see next that the subspaces $\H_{(\E,\F)}(\SS)$ and $\F(\SS)$ actually provide an orthogonal decomposition of the domain of the resistance forms $(\E,\F)\in\RF_S$.

\thm\label{IRF.thm10}
Let $(\E, \F) \in \RF_S$ and let $\{(\E_m, \F_m, \eta_m)\}_{m \ge 0} \subseteq \RF_S^{(0)} \times (0, \infty)$ be the resolution of $(\E, \F)$. Then,
\begin{equation}\label{IRF.eq70}
\F = \H_{(\E, \F)}(\SS) \oplus \F(\SS)
\end{equation}
and for any $u \in \F$,
\begin{equation}\label{IRF.eq80}
\E(u, u) = \E|_{\SS}(u|_{\SS}, u|_{\SS}) + \sum_{k = 1}^{\infty}\frac 1{\eta_k}\D_k^I(u - h_{\SS}(u), u - h_{\SS}(u)),
\end{equation}
where $h_{\SS}(u)$ denotes the $\SS$-harmonic extension of $u|_\SS$.
\endthm

\demo
The direct sum decomposition~\eqref{IRF.eq70} follows from Proposition~\ref{BRF.prop30}. Combining Proposition~\ref{BRF.prop30} and Lemma~\ref{IRF.lemma50} we immediately verify \eqref{IRF.eq80}.
\enddemo

\cor\label{IRF.cor10}
Let $(\E, \F),(\E', \F') \in \RF_S$. Then, the following conditions are equivalent:
\begin{itemize}[leftmargin=.325in]
\item[\rm{(E1)}] $(\E, \F) = (\E', \F')$.
\item[\rm{(E2)}] $(\E|_{V_m}, \l(V_m)) = (\E'|_{V_m}, \l(V_m))$ for any $m \ge 0$.
\item[\rm{(E3)}] $(\E|_{\SS}, \F|_{\SS}) = (\E'|_{\SS}, \F'|_{\SS})$.
\item[\rm{(E4)}] $\R_{(\E, \F)} = \R_{(\E', \F')}$ and $r_0(\E, \F) = r_0(\E', \F')$.
\end{itemize}
\endcor

\demo
By Theorem~\ref{BRF.thm20}, (E2) implies (E3). Since $(\E|_{\SS})|_{V_m} = \E|_{V_m}$ and $(\E'|_{\SS})|_{V_m} = \E'|_{V_m}$, we see that (E3) implies (E2). In view of Lemma~\ref{IRF.lemma30}, (E4) implies (E2). Conversely, assume that (E2) holds. Since $\E|_{V_0} = \E'|_{V_0}$, then $r_0(\E, \F) = r_0(\E', \F')$ and by Lemma~\ref{IRF.lemma30} we obtain $E_{\R_{(\E, \F)}, m} = E_{\R_{(\E', \F')}, m}$ for any $m \ge 1$. Therefore, it follows that $\R_{(\E, \F)} = \R_{(\E', \F')}$ and hence we have (E4). Moreover, (E1) immediately implies (E3) and it only remains to verify that (E3) implies (E1). \\
Let us assume (E2), (E3) and (E4). By Lemma~\ref{IRF.lemma40}, if $h_{\SS}\colon\F|_{\SS} \to \F$ and $h'|_{\SS}\colon\F'|_{\SS} \to \F'$ are the $\SS$-harmonic extension maps associated with $(\E, \F)$ and $(\E', \F')$ respectively, then $h_{\SS} = h'_{\SS}$ and hence $\H_{(\E, \F)}(\SS) = \H_{(\E', \F')}(\SS)$. If $\{(\E_m, \F_m, \eta_m)\}_{m \ge 0}$ and $\{(\E'_m, \F'_m, \eta'_m)\}_{m \ge 0}$ are the resolutions of $(\E, \F)$ and $(\E', \F')$ respectively, then~\eqref{IRF.eq40} implies that $\eta_m = \eta'_m$. Lemma~\ref{IRF.lemma50} thus yields $\F(\SS) = \F'(\SS)$ and by~\eqref{IRF.eq70} it follows that $\F = \F'$. Finally, since $(\E|_{\SS}, \F|_{\SS}) = (\E'|_{\SS}, \F'|_{\SS})$ and $\eta_m = \eta'_m$ for any $m \ge 0$, \eqref{IRF.eq80} shows that $\E(u, u) = \E'(u, u)$ for any $u \in \F = \F'$, hence $(\E, \F) = (\E', \F')$.
\enddemo

Resistance forms on $K$ constructed by means of sequences of matching pairs as explained in Section~\ref{CRF} have the property of belonging to $\RF_S^N$.

\lemma\label{IRF.lemma60}
Let $\R \in \MP^{\BbN}$. Then $r_0(\E_{\R}, \F_{\R}) = 1$.
\endlemma

\demo
Let $f: V_0 \to \BbR$ and let $u$ be the $V_0$-harmonic function  with respect to $(\hE_{\R}, \hF_{\R})$ with boundary value $f$. Since $(\hE_{\R}|_{V_0}, \hF_{\R}|_{V_0}) = (Q_0^{\SS}, \l(V_0))$, we have
\[
\hE_{\R}(u, u) = Q_0^{\SS}(f, f).
\]
By Lemma~\ref{RNS.lemma100}, $\hF_{\R} \subseteq C(\SS)$ and hence $u \in C(\SS)$. For any $(w, (i, j)) \in W_* \times B$, we extend the domain of $u$ to each $e_{ij}^w$ by defining
\[
\vp|_{e_{ij}^w}((1 - t)G_w(p_{ij}) + tG_w(p_{ji})) = (1 - t)u(G_w(p_{ij})) + tu(G_w(p_{ji}))
\]
for any $t\in [0, 1]$ and $\vp|_{\SS} =u$. In this manner, $u$ is extended to a function $\vp$ on $K$ such that $\vp \in \tF$. Since for any $m\geq 1$
\[
Q_0^{\SS}(f, f) = E_{\R, m}(\vp, \vp) = \E_{\R, m}(\vp, \vp),
\]
$\vp \in \F_{\R}$ and $\E_{\R}(\vp, \vp) = Q_0^{\SS}(f, f)$. Now, if  $v \in \F_{\R}$ and $v|_{V_0} = f$,~\eqref{RNS.eq105} yields
\[
\E_{\R}(\vp, \vp) = Q_0^{\SS}(f, f) \le E_{\R, m}(v,v) \le \E_{\R, m}(v,v) \le \E_{\R}(v,v).
\]
Therefore, $\vp$ is the $V_0$-harmonic function with respect to $(\E_{\R}, \F_{\R})$ with boundary value $f$ and hence $\E_{\R}|_{V_0} = Q_0^{\SS}$. This shows $r_0(\E_{\R}, \F_{\R}) = 1$.
\enddemo

We conclude this section by showing that $\RF_S = \{(\delta\E_{\R}, \F_{\R})~|~\R \in \MP^{\BbN}, \delta > 0\}$.
\thm\label{IRF.thm20}
\begin{itemize}[leftmargin=0.25in]
\item[\rm (1)] $\R_{(\E_\R, \F_\R)} = \R$ for any $\R \in \MP^{\BbN}$.
\vspace*{-.1in}
\item[\rm (2)] $(\E_{\R_{(\E, \F)}}, \F_{\R_{(\E, \F)}}) = (r_0(\E, \F)\E, \F)$ for any $(\E, \F) \in \RF_S$.
\end{itemize}
\noindent
In particular, the map $\R \in \MP^{\BbN} \to (\E_{\R}, \F_{\R}) \in \RF_S^N$ is bijective.
\endthm

\demo
(1)~Let $\R = \{(r_m, \rho_m)\}_{m \ge 1}$. By Theorem~\ref{thm20}, the resolution of $(\E_{\R}, \F_{\R})$ is given by $\{((\delta_m)^{-1}\E_{\R^{(m)}}, \F_{\R^{(m)}}, \c_m)\}_{m \ge 0}$, where $\delta_m = r_1{\cdots}r_m$ and $\c_m = \delta_{m - 1}\rho_m$. By Lemma~\ref{IRF.lemma60}, $r_0((\delta_m)^{-1}\E_{\R^{(m)}}, \F_{\R^{(m)}}) = \delta_m$ and therefore $\R_{(\E_{\R}, \F_{\R})} = \R$.\\
(2)~Without loss of generality, we may assume that $r_0(\E, \F) = 1$. Set $\R = \R_{(\E, \F)}$. Applying (1), $\R_{(\E_{\R}, \F_{\R})} = \R$ and by Lemma~\ref{IRF.lemma60}, $r_0(\E_{\R}, \F_{\R}) = 1 = r_0(\E, \F)$. Thus, condition (E4) in Corollary~\ref{IRF.cor10} is satisfied and hence $(\E_{\R}, \F_{\R}) = (\E, \F)$.
\enddemo

\remark
This result reveals that $\RF_S^N$ is in fact the set of fixed points of the mapping
\begin{align*}
\Phi\colon&\;\RF_S\;\longrightarrow\;\RF_S\\
&(\E,\F)\;\mapsto (r_0(\E,\F)^{-1}\E_{\R_{(\E,\F)}},\F_{\R_{(\E,\F)}}).
\end{align*}
\endremark

\setcounter{equation}{0}

\section{Classification of resistance forms derived from \\ matching pairs}\label{CMP}

In the previous section we have identified any complete symmetric resistance form on SSG with a resistance form $(\E_\R,\F_\R)$ derived from a sequence of matching pairs $\R$ up to multiplication by a constant. The present section analyzes the detailed structure of $(\E_\R,\F_\R)$ by decomposing it into two parts, called the \textit{SG part} in allusion to the reminiscence of SG in SSG, and the \textit{line part} that corresponds to the cable system/quantum graph approach. In Theorem~\ref{CRF.thm10} we use a certain property of the sequence $\R$ to determine when the SG part is non-trivial and therefore captures the reminiscence of the SG in the geometric structure of SSG.\par
Let $\R = \{(r_m, \rho_m)\}_{m \ge 1} \in \MP^{\BbN}$ and set $\delta_m = r_1{\cdots}r_m$ and $\c_m = \delta_{m - 1}\rho_m$. For any $u \in \F_{\R}$ and $m \ge 1$,
\[
\sum_{k = 1}^m \frac 1{\c_k}\D_k^I(u,u) \le \E_{\R}(u, u)
\]
and the left-hand side is monotonically increasing with respect to $m$. We can thus define $\E^I_{\R}(u, u)$ as 
\begin{equation}\label{IND.eq100}
\E_{\R}^I(u, u) = \sum_{k = 1}^{\infty} \frac 1{\c_k}\D_k^I(u,u).
\end{equation}
Moreover, we define
\begin{equation}\label{IND.eq200}
\E_{\R}^{\SS}(u, u) = \E_{\R}(u, u) - \E^I_{\R}(u, u)
\end{equation}
and therefore 
\[
\E_{\R}^{\SS}(u, u) = \lim_{m \to \infty} \frac 1{\delta_m}Q_m^{\SS}(u, u).
\]
Notice that both $\E_{\R}^{\SS}$ and $\E_{\R}^I$ are non-negative quadratic forms on $\F_{\R}$. As a consequence, it follows that
\[
\E_{\R}(u, v) = \E_{\R}^{\SS}(u, v) + \E_{\R}^I(u, v).
\]
We call $\E^{\SS}_{\R}$ (resp. $\E^I_{\R}$) the SG part (resp. the line part) of $\E_{\R}$.
It is easy to see that, in order for $\E_\R$ to be a resistance form, the line part should be non-zero. On the contrary, the SG part may vanish, as we will see in the the course of our discussion.\par
The following useful lemma is an exercise of undergraduate calculus. 
\lemma\label{lemma30}
Let $\{(r_m, \rho_m)\}_{m \ge 1} \in \MP^{\BbN}$ and define
\[
\kappa_m = -\log{(1 - \rho_m)}.
\]
The following three conditions are equivalent:
\begin{itemize}[leftmargin=.25in]
\item[\rm{(a)}] $\sum_{m = 1}^{\infty} \kappa_m < +\infty$,
\item[\rm{(b)}] $\sum_{m = 1}^{\infty} \rho_m < +\infty$,
\item[\rm{(c)}] There exists $C > 0$ such that
\[
C\Big(\frac35\Big)^m \le r_1r_2\cdots{r_m}
\]
for any $m \ge 1$.
\end{itemize}
\endlemma

Note that $\kappa_m = \log{\frac 35} - \log{r_m}$ for each $m\geq 1$. Thus if $\sum_{m = 1}^{\infty} 
\rho_m < +\infty$ and we set $\kappa = \sum_{m = 1}^{\infty} \kappa_m$, then
\[
r_1r_2{\cdots}r_m = e^{-\kappa}e^{\sum_{i \ge m + 1} \kappa_i}\Big(\frac 35\Big)^m.
\]
With this equality one is led to the following lemma.

\lemma\label{lemma40}
Let $\{(r_m, \rho_m)\}_{m \ge 1} \in \MP^{\BbN}$. There exist $C_* > 0$ and a sequence $\{c_m\}_{m \ge 1}$ such that $\lim_{m \to \infty} c_m = 0$ and
\[
\frac 1{r_1r_2\cdots{r_m}} = C_*\Big(\frac 53\Big)^m(1 - c_m)
\]
if and only if $\sum_{m = 1}^{\infty} \rho_m < \infty$. Moreover, if $\sum_{m = 1}^{\infty} \rho_m < \infty$, then $C_* = \prod_{m = 1}^\infty (1 - \rho_m)^{-1}$, $c_m \ge 0$ for any $m \ge 0$, and $\{c_m\}_{m \ge 1}$ is monotonically decreasing.
\endlemma

\definition\label{CRF.def00}
Define $\F_{\R, *} = \F_{\R} \cap C(K_*)$.
\enddefinition

As announced at the beginning of the section, the next theorem reveals the sufficient condition for the survival of $\E_{\R}^{\SS}$, that is $\sum_{m = 1}^{\infty} \rho_m < \infty$. In the next section, this condition will be seen to be necessary as well.

\thm\label{CRF.thm10}
Let $\R = \{(r_m, \rho_m)\}_{m\geq 1}\in \MP^{\BbN}$. $\F_{\R, *}$ consists of constants if and only if $\sum_{m=1}^\infty \rho_m = \infty$. Furthermore, if $\sum_{m=1}^\infty \rho_m < \infty$, then $\F^* \subseteq \F_{\R, *}$ and
\[
\E_{\R}(u, u) = \E_{\R}^{\SS}(u, u) = C_*\E^*(u, u)
\]
for any $u \in \F^*$, where $C_* = \prod_{m=1}^\infty(1 - \rho_m)^{-1}$.
\endthm

\demo
For each $m\geq 1$, set $\delta_m = r_1{\cdots}r_m$. If $\sum_{m = 1}^{\infty} \rho_m < \infty$, then Lemma~\ref{lemma40} yields
\[
\E_{\R, m}(u, u) = \frac 1{\delta_m}Q^{\SS}_m(u, u) = C_*(1 - c_m)\Big(\frac 53\Big)^mQ^{\SS}_m(u, u)
\]
for any $u \in \F^*$. By Theorem~\ref{SG.thm10} we now have that
\[
\E_{\R}(u, u) = \lim_{m \to \infty} \frac 1{\delta_m}Q^{\SS}_m(u, u) = C_*\E^*(u, u)
\]
and hence $\F^* \subseteq \F_{\R, *}$. On the other hand, if $\sum_{m=1}^\infty \rho_m = \infty$, then for any $u\in C(K_*)$, 
\[
\E_{\R, m}(u, u) = \frac 1{\delta_m}Q^{\SS}_m(u, u) = \frac 1{\delta_m}\Big(\frac 35\Big)^m\cdot\Big(\frac 53\Big)^mQ^{\SS}_m(u, u).
\]
From Lemma~\ref{lemma30}, it follows that $\frac 1{\delta_m}\big(\frac 35\big)^m$ is unbounded as $m \to \infty$ and by Proposition~\ref{SG.prop10}, $\lim_{m \to \infty}\big(\frac 53\big)^mQ^{\SS}_m(u, u) > 0$ unless $u$ is constant. Therefore, $u \in \F_{\R}$ if and only if $u$ is constant.
\enddemo

We conclude this section with some preliminary results concerning the domains of the SG part and the line part.

\definition\label{CMP.def10}
Let $\R \in \MP^{\BbN}$. Define
\begin{align*}
\F_{\R}^I &= \{u~|~u \in \F_{\R},~\E_{\R}^{\SS}(u, u) = 0\},\\
\F_{\R}^{\SS} &= \{u~|~u \in \F_{\R},~\E_{\R}^I(u, u) = 0\}.
\end{align*}
\enddefinition

\lemma\label{CMP.lemma10}
Let $\R \in \MP^{\BbN}$. Then,
\[
\F_{\R}^I = \{u~|~u \in \F_{\R},~\E_{\R}^{\SS}(u, v) = 0\text{ for any }v \in \F_{\R}\}
\]
and
\[
\F_{\R}^{\SS} = \{u~|~u \in \F_{\R},~\E_{\R}^I(u, v) = 0\text{ for any }v \in \F_{\R}\}.
\]
\endlemma

\demo
Let us consider first $\F_{\R}^I$. Applying the Cauchy-Schwartz inequality,
\[
\E_{\R}^{\SS}(u, v)^2 \le \E_{\R}^{\SS}(u, u)\E_{\R}^{\SS}(v, v)
\]
for any $u,v\in\F_\R$. Therefore, if $u \in \F_{\R}^I$, then $\E_{\R}^{\SS}(u, v) = 0$ for any $v \in \F_\R$. The converse direction is obvious. The argument for $\F_{\R}^{\SS}$ works verbatim.
\enddemo

\lemma\label{CMP.lemma15}
Let $\R \in \MP^{\BbN}$. Then, $\F_{\R}^{\SS} = \F_{\R} \cap C(K_*)$.
\endlemma

\demo
If $u \in \F_{\R} \cap C(K_*)$, then $u$ is constant on every $e_{ij}^w$. Therefore, $\E_{\R}^I(u, u) = 0$ and hence $u \in \F_{\R}^{\SS}$. Conversely, assume that $u \in \F_{\R}^{\SS}$. Then, $\D_{e_{ij}^w}(u|_{e_{ij}^w}, u|_{e_{ij}^w}) = 0$ for every $e_{ij}^w$ and hence $u$ is constant on $e_{ij}^w$. Thus, $u \in C(K_*)$.
\enddemo

\setcounter{equation}{0}

\section{Projection to the line part}\label{SRF}

Let us define $(\MP^{\BbN})^I \subseteq \MP^{\BbN}$ by
 \begin{equation}\label{IND.eq300}
(\MP^{\BbN})^I = \{\R~|~\R \in \MP^{\BbN},~(\E_{\R}, \F_{\R}) = (\E_{\R}^I, \F_{\R}^I)\}.
\end{equation}
In this section, we are going to introduce a natural projection $\L\colon\MP^{\BbN} \to (\MP^{\BbN})^I$ and through an explicit expression of this mapping, it will be shown in Theorem~\ref{SRF.thm30} that
\[
(\MP^{\BbN})^I = \big\{\R~\big|~ \R = \{(r_m, \rho_m)\}_{m \ge 1} \in \MP^{\BbN}, \sum_{m = 1}^{\infty} \rho_m = \infty\big\}.
\]
In other words, the converse of Theorem~\ref{CRF.thm10} holds, i.e. $\sum_{m=1}^\infty\rho_m=\infty$ if and only if $\E_\R=\E_\R^I$, or equivalently $\E_\R^\SS=0$.

Before doing so, we present some results concerning a general theory that is not confined to resistance forms on $K$ and is applicable in very abstract settings.

Let $(\E, \F)$ be a resistance form on a set $X$ and let $R$ be the associated resistance metric. Assume that there exist non-negative symmetric quadratic forms $\E^{(1)}(\cdot, \cdot)$ and  $\E^{(2)}(\cdot, \cdot)$ on $\F \times \F$ such that
\[
\E(u, v) = \E^{(1)}(u, v) + \E^{(2)}(u, v)
\]
for any $u, v \in \F$. Define $\F^{(2)} = \{u~|~u \in \F,~\E^{(1)}(u, u) = 0\}$ and note that
\[
\E(u, v) = \E^{(2)}(u, v)
\]
for any $u, v \in \F^{(2)}$. As in Definition~\ref{BRF.def10}, for any $u, v \in \F$ we define $u \sim v$ if and only if $u-v$ is constant.

\lemma\label{SRF.lemma10}
$(\F^{(2)}/\!\!\sim, \E)$ is a closed subspace of $(\F/\!\!\sim, \E)$.
\endlemma

\demo
Let $x \in X$ and set $\F^{(2)}_x = \{u|u \in \F^{(2)}, u(x) = 0\}$ and $\F_x = \{u|u \in \F, u(x) = 0\}$. It suffices to show that $(\F^{(2)}_x, \E)$ is a closed subspace of $(\F_x, \E)$. Let $\{u_n\}_{n \ge 1} \in \F^{(2)}_x$ and suppose that $\E(u_n - u, u_n - u) \to 0$ as $n \to \infty$ for some $u \in \F_x$. Then,
\[
\E^{(1)}(u - u_n, u - u_n) = \E^{(1)}(u_n, u_n) - 2\E^{(1)}(u_n, u) + \E^{(1)}(u, u) \to 0
\]
as $n \to \infty$. Since $\E^{(1)}(u_n, u_n) = 0$ and $\E^{(1)}(u_n, u)^2 \le \E^{(1)}(u_n, u_n)\E^{(1)}(u, u)$, it follows that $\E^{(1)}(u, u) = 0$.  Thus, $(\F^{(2)}_x, \E)$ is closed and so is $(\F^{(2)}/\!\!\sim, \E)$.
\enddemo
Using the above lemma, we may easily verify the following statement.
\thm\label{SRF.thm10}
$(\E, \F^{(2)})$ is a resistance form on $X$ if the following two conditions are satisfied:
\begin{itemize}[leftmargin=0.25in]
\item[\rm{(1)}] For any $x \neq y \in X$, there exists $u \in \F^{(2)}$ such that $u(x) \neq u(y)$.
 \item[\rm{(2)}] For any $u \in \F$, $\E^{(1)}(\overline{u}, \overline{u}) \le \E^{(1)}(u, u)$ and $\E^{(2)}(\overline{u}, \overline{u}) \le \E^{(2)}(u, u)$.
\end{itemize}
\endthm

\begin{proof}
(RF1) and (RF4) hold because $\F^{(2)}\subseteq\F$ and $(\E,\F)$ is a resistance form. (RF3) is condition (1) and (RF5) is condition (2). To prove (RF2), notice that $(\F/\!\!\sim,\E)$ is complete because $(\E,\F)$ is a resistance form. By Lemma~\ref{SRF.lemma10}, $(\F^{(2)}/\!\!\sim,\E)$ is a closed subspace of $(\F/\!\!\sim,\E)$ and hence also complete. 
\end{proof}

 Back to resistance forms on SSG, we start by constructing the projection $\L$ from $\MP^\BbN$ onto $(\MP^\BbN)^I$.

\thm\label{SRF.thm20}
Let $\R = \{(r_m, \rho_m)\}_{m \ge 1} \in \MP^{\BbN}$.  If 
\[
\F_{\R}^I = \{u~|~u \in \F_{\R}, \E_{\R}^{\SS}(u, u) = 0\},
\]
then $(\E_{\R}^I, \F_{\R}^I) \in \RF_S$ and its resolution is given by $\{((\delta_m)^{-1}\E_{\R^{(m)}}^I, \F_{\R^{(m)}}^I, \c_m)\}_{m \ge 0}$, where $\delta_m = r_1{\cdots}r_m$ and $\c_m = \delta_{m - 1}\rho_m$ for any $m \ge 1$. Moreover, $r_0(\E_{\R}^I, \F_{\R}^I) \le 1$ and there exists a unique $\R' \in \MP^{\BbN}$ such that $(r_0(\E_{\R}^I, \F_{\R}^I)\E_{\R}^I, \F_{\R}^I) = (\E_{\R'}, \F_{\R'})$.
\endthm
To prove this theorem, we need the following lemma.

\lemma\label{SRF.lemma15}
\[
\F_{\R}^I=\{u~|~u\in C(K),~u{\circ}G_i \in \F_{\R^{(1)}}^I\text{ for any }i \in S,~u|_{e_{ij}} \in H^1(e_{ij}) \text{ for any }(i, j) \in B\}
\]
and
\begin{equation}\label{SRF.eq50}
\E_{\R}^I(u, u) = \sum_{i \in S} \frac 1{r_1}\E_{\R^{(1)}}^I(u{\circ}G_i, u{\circ}G_i) + \frac{1}{\rho_1}\D_1^I(u, u)
\end{equation}
for any $u \in \F_{\R}^I$.
\endlemma

\demo
Note that for any $m \ge 1$ and $u\in\F_\R^I$,
\[
\frac 1{\delta_m}Q_m^{\SS}(u, u) = \frac 1{\delta_m}\sum_{i \in S}Q_{m - 1}^{\SS}(u{\circ}G_i, u{\circ}G_i).
\]
By definition, $u \in \F_{\R}^I$ if and only if $u \in \F_{\R}$ and
\[
\lim_{m \to \infty} \frac 1{\delta_{m + 1}}Q_m^{\SS}(u{\circ}G_i, u{\circ}G_i) = 0
\]
for any $i \in S$, i.e. $u \in \F_{\R}^I$ if and only if $u \in \F_{\R}$ and $u{\circ}G_i \in \F_{\R^{(1)}}^I$. This immediately implies the desired equivalence. 
Since $\{((\delta_m)^{-1}\E_{\R^{(m)}}, \F_{\R^{(m)}}, \c_m)\}_{m \ge 0}$ is the resolution of $(\E_{\R}, \F_{\R})$, we obtain \eqref{SRF.eq50}.
\enddemo
\demo[Proof of Theorem~\ref{SRF.thm20}]
Applying Theorem~\ref{SRF.thm10} to $(\E_{\R}, \F_{\R})$ with $\E^{(1)} = \E_{\R}^{\SS}$ and $\E^{(2)} = \E_{\R}^I$, we see that $(\E_{\R}^I, \F_{\R}^I)$ is a resistance form on $K$.\\
Moreover, by Lemma~\ref{SRF.lemma15}, $u \in \F_{\R}^I$ if and only if $u|_{G_i(K)} \in \F_{\R}^I|_{G_i(K)}$ for any $i \in S$ and $u|_{e_{ij}} \in H^1(e_{ij})$. 
If $R(\cdot, \cdot)$ and $R^I(\cdot, \cdot)$ are the resistance metrics associated with $(\E_{\R}, \F_{\R})$ and $(\E_{\R}^I, \F_{\R}^I)$ respectively, then
\begin{equation}\label{SRF.eq90}
R^I(x, y) = \sup\Big\{\frac{|u(x) - u(y)|^2}{\E_{\R}(u, u)}~\big|~u \in \F_{\R}^I, u(x) \neq u(y)\Big\} \le R(x, y)
\end{equation}
for any $x, y \in K$. Hence, the identity map $\iota$ from $(K, R)$ to $(K, R^I)$ is continuous and since $(K, R)$ is compact, the map $\iota$ is a homeomorphism. Furthermore, $\E_\R^I$ is invariant under all geometric symmetries of $K$ because $\E_\R$ is. Combining these previous facts, we conclude that $(\E_{\R}^I, \F_{\R}^I) \in \RF_S^{(0)}$. Applying~\eqref{SRF.eq50} to $(\E_{\R^{(m)}}^I, \F_{\R^{(m)}}^I)$ repeatedly, we get that $(\E_{\R}^I, \F_{\R}^I) \in \RF_S$ and its resolution is $\{((\delta_m)^{-1}\E_{\R^{(m)}}^I, \F_{\R^{(m)}}^I, \c_m)\}_{m \ge 0}$. Finally,
\[
r_0(\E_{\R}^I, \F_{\R}^I) = \frac 32R^I(p_1, p_2) \le \frac 32R(p_1, p_2) = 1
\]
and the existence of $\R' \in \MP^{\BbN}$ follows immediately from Theorem~\ref{IRF.thm20}.
\enddemo

\definition\label{SRF.def10}
For any $\R \in \MP^{\BbN}$, $\L(\R) \in \MP^{\BbN}$ is defined as $\R'$ given in Theorem~\ref{SRF.thm20}.
\enddefinition

\lemma\label{SRF.lemma17}
Let $\R \in \MP^{\BbN}$. If $\L(\R) = \{(s_m, \s_m)\}_{m \ge 1}$, then 
\begin{equation}\label{SRF.eq160}
\sum_{m = 1}^{\infty} \s_m = \infty.
\end{equation}
\endlemma

\demo
Notice that for any $u \in \F^* \cap \F_{\L(\R)}$, $u$ is constant on every $e^w_{ij}$ and hence $\E_{\L(\R)}(u, u) = r_0(\E_{\R}^I, \F_{\R}^I)\E_{\R}^I(u, u) = 0$. Since $(\E_{\L(\R)}, \F_{\L(\R)})$ is a resistance form, $u$ is constant on $K$. By Theorem~\ref{CRF.thm10}, we see that $\sum_{m = 1}^{\infty} \s_m = \infty$.
\enddemo

The next lemma gives an explicit expression of $\L(\R)$, which plays an essential role in the rest of the section.

\lemma\label{SRF.lemma20}
Let $\R = \{(r_m, \rho_m)\}_{m\ge 1} \in \MP^{\BbN}$ and let $\rho_0 = r_0(\E_{\R}^I, \F_{\R}^I)$. If $\L(\R)= \{(s_m, \s_m)\}_{m \ge 1}$, then
\begin{equation}\label{SRF.eq120}
\rho_0\prod_{i = 1}^{m - 1} (1 - \s_i) = \prod_{i = 1}^{m - 1} (1 - \rho_i) - (1 - \rho_0)
\end{equation}
for any $m \ge 1$. In particular,
\begin{equation}\label{SRF.eq130}
\s_m = \frac{\prod_{i = 1}^{m - 1}(1 - \rho_i)}{\prod_{i = 1}^{m - 1} (1 - \rho_i) - (1 - \rho_0)}\rho_m
\end{equation}
and
\begin{equation}\label{SRF.eq135}
\rho_m = \frac{\rho_0\prod_{i = 1}^{m - 1}(1 - \s_i)}{\rho_0\prod_{i = 1}^{m - 1} (1 - \s_i) + (1 - \rho_0)}\s_m.
\end{equation}
\endlemma

\demo
For $(w,(i, j)) \in W_*\times B$, choose $u \in \F_{\R}$ so that $u(x) = 0$ for any $x \notin e^w_{ij}$ and $\E_{\R}(u, u) > 0$. Then, $u \in \F_{\R}^I$ and $\rho_0\E_{\R}^I(u, u) = \E_{\R'}(u, u) = \E_{\R'}^I(u, u)$. Hence, we get
\begin{equation}\label{SRF.eq100}
r_1r_2\cdots{r_{m - 1}}\rho_m = \rho_0s_1s_2\ldots{s_{m - 1}}\s_m
\end{equation}
for any $m \ge 1$, which yields
\begin{equation}\label{SRF.eq110}
\rho_m\prod_{i = 1}^{m - 1}(1 - \rho_i) = \rho_0\s_m\prod_{i = 1}^{m- 1}(1 - \s_i).
\end{equation}
By induction, we obtain \eqref{SRF.eq120}.
\enddemo

\lemma\label{SRF.lemma30}
Let $\R = \{(r_m, \rho_m)\}_{m\ge 1} \in \MP^{\BbN}$ and let $\rho_0 = r_0(\E_{\R}^I, \F_{\R}^I)$.  If $\L(\R) = \{(s_m, \s_m)\}_{m\ge 1}$, then
\begin{equation}\label{SRF.eq140}
\s_m = \frac{1 - \a_m}{\rho_0 - \a_m}\rho_m
\end{equation}
for any $m\ge 1$, where $\a_m = 1 - \prod_{i = 1}^{m- 1}(1 - \rho_i)$. In particular, 
\[
\rho_0 \ge \lim_{m \to \infty} \a_m.
\]
\endlemma

\begin{proof}
The equality follows directly from Lemma~\ref{SRF.lemma20}, which also implies that for any $m\geq 1$,
\[
\rho_0=\rho_0\prod_{i=1}^{m-1}(1-\sigma_i)+1-\prod_{i=1}^{m-1}(1-\rho_i)\geq 1-\prod_{i=1}
^{m-1}(1-\rho_i)=\alpha_m
\]
and therefore $\rho_0\geq\lim_{m\to\infty}\alpha_m$.
\end{proof}

\remark
$\{\a_n\}_{n \ge 1}$ is monotonically increasing and $\a_n \uparrow \a$ as $n \to \infty$ for some $\a \in (0, 1]$.
\endremark 

Finally, we present the main theorem of this section. It characterizes $(\MP^\BbN)^I$ and essentially says that the SG part $\E_{\R}^{\SS}$ truly exists if and only if $\sum_{m = 1}^\infty \rho_m < \infty$.

\thm\label{SRF.thm30}
Let $\R = \{(r_m, \rho_m)\}_{m \ge 1} \in \MP^{\BbN}$. Then, $\L(\R) = \R$ if and only if $\sum_{m = 1}^\infty \rho_m = \infty$. In particular, $\E_{\R} = \E_{\R}^I$ if and only if $\sum_{m =1}^\infty \rho_m = \infty$.
\endthm

\demo
Assume that $\sum_{m = 1}^{\infty} \rho_m = \infty$. Then $\a = 1$, which implies that $\rho_0 \ge 1$ and therefore $\rho_0 = 1$. In view of~\eqref{SRF.eq130}, we have that $\rho_m = \s_m$ for any $m \ge 1$, hence $\R = \R'$. Thus we have shown that $\E_{\R} = \E_{\R'} = \E_{\R}^I$. Conversely, if $\R = \R'$, then Lemma~\ref{SRF.lemma17} shows that $\sum_{m = 1}^{\infty} \rho_m= \sum_{m = 1}^{\infty} \s_m = \infty$.
\enddemo

As a consequence of this theorem,
\[
\L(\L(\R)) = \L(\R)
\]
for any $\R \in \MP^{\BbN}$ and $\L(\MP^{\BbN}) = (\MP^{\BbN})^I$. Thus, we may regard $\L$ as a projection onto $(\MP^{\BbN})^I$.

We finish this section with several useful equalities leading to an explicit expression of $r_0(\E_\R^I,\F_\R^I)$ in terms of the elements of $\R$.

\lemma\label{SRF.lemma40}
Let $\R=\{(r_m, \rho_m)\}_{m \ge 1} \in \MP^{\BbN}$. Then,
\begin{itemize}[leftmargin=.25in]
\item[\rm{(1)}] 
For any $m\geq 1$,
\[
\sum_{i = 1}^{m} \bigg(\frac 53\bigg)^{i - 1}\c_i + \bigg(\frac 53\bigg)^m\delta_m = 1,
\] 
where $\delta_m = r_1{\cdots}r_m$ and $\c_i = \delta_{i - 1}\rho_i$.
\item[\rm{(2)}] $\sum_{m = 1}^\infty \rho_m = \infty$ if and only if $\lim_{m \to \infty}\big(\frac 53\big)^m\delta_m = 0$.
\end{itemize}
\endlemma

\demo
Let $\c_m = r_1{\cdots}r_{m - 1}\rho_m$. Since $\frac 53r_m + \rho_m = 1$, we have
\[\begin{split}
\Big(\frac 53\Big)^{m - 1}\c_m &= (1 - \rho_1){\cdots}(1 - \rho_{m - 1})\rho_m\\
& = (1 - \rho_1){\cdots}(1 - \rho_{m - 1}) - (1 - \rho_1){\cdots}(1 - \rho_m)
\end{split}\]
and hence
\[
\sum_{i = 1}^m\bigg(\frac 53\bigg)^{i - 1}\c_i = 1 - \prod_{i = 1}^m (1 - \rho_i)= 1 - \bigg(\frac 53\bigg)^{m}\delta_m.
\]
This proves (1). Assertion (2) follows immediately from the fact that $\prod\limits_{i = 1}^m (1 - \rho_i) = \big(\frac 53\big)^{m}\delta_m$.
\enddemo

\prop\label{SRF.prop10}
Let $\R = \{(r_m, \rho_m)\}_{m\ge 1} \in \MP^{\BbN}$. Then,
\begin{equation}\label{SRF.eq150}
r_0(\E_{\R}^I, \F_{\R}^I) = 1 - \prod_{m = 1}^{\infty} (1 - \rho_m) = \sum_{m = 1}^{\infty} \bigg(\rho_m\prod_{i = 1}^{m - 1} (1 - \rho_i)\bigg).
\end{equation}
In particular, $r_0(\E_{\R}^I, \F_{\R}^I) < 1$ if and only if $\sum_{m = 1}^{\infty}\rho_m < \infty$.
\endprop
 
\demo
Set $\rho_0 = r_0(\E_{\R}^I, \F_{\R}^I)$ and $\L(\R) = \{(s_m, \s_m)\}_{m \ge 1}$. If $\sum_{m = 1}^{\infty} \rho_m = \infty$, then we have already shown in the proof of Theorem~\ref{SRF.thm30} that $\rho_0 = 1$. Since $\prod_{m = 1}^{\infty}(1- \rho_m) = 0$, Lemma~\ref{SRF.lemma40} implies~\eqref{SRF.eq150}.\par
Suppose that $\sum_{m = 1}^{\infty} \rho_m < \infty$ and set $\a = 1 - \prod_{m = 1}^{\infty} (1 - \rho_m)$. Note that $(1 - \a_m)\rho_m = (5/3)^{m - 1}\c_m$ as in the proof of Lemma~\ref{SRF.lemma40}-(1).  Therefore, if $\rho_0  > \a$, then~\eqref{SRF.eq140} and Lemma~\ref{SRF.lemma40}-(1) lead to
\[
\sum_{m = 1}^{\infty} \s_m \le \sum_{m  = 1}^{\infty}\frac{(1 - \a_m)\rho_m}{\rho_0 - \a} \le \frac 1{\rho_0 - \a} < \infty.
\]
This contradicts~\eqref{SRF.eq160}, hence $\rho_0 = \a$. Applying Lemma~\ref{SRF.lemma40} again, we immediately obtain \eqref{SRF.eq150}.
\enddemo
\setcounter{equation}{0}

\section{Domain of resistance forms given by infinite \\ sequences of matching pairs}\label{DRF}

 The results obtained in previous sections come together in the present one to prove the main theorem of this paper, Theorem~\ref{RFF.thm00}. In fact, Theorem~\ref{IRF.thm20} and Theorem~\ref{SRF.thm30} already identify any completely symmetric resistance form on SSG as the sum of its line part and its SG part, whenever the latter survives. This identification is now completed by giving a full description of the domains of these forms. This characterization of the domains in the next theorem is the key step to showing Theorem~\ref{RFF.thm00}.

\thm\label{DRF.thm10}
Let $\R = \{(r_m, \rho_m)\}_{m \ge 1} \in \MP^{\BbN}$ and set $R_* = \prod_{m = 1}^{\infty} (1 - \rho_m)$. Moreover, define
\[
\eta_m = \frac{r_1{\cdots}r_{m - 1}\rho_m}{1 -R_*}
\]
for any $m \ge 1$ and $\eta = \{\eta_m\}_{m \ge 1}$.
\begin{itemize}[leftmargin=.25in]
\item[\rm(1)] If $\sum_{m = 1}^{\infty} \rho_m = \infty$, then $R_* = 0$, $\F_{\R} = \F_{\eta}$ and 
\[
\E_{\R}(u, v) = \D_{\eta}^I(u, v)
\]
for any $u, v \in \F_{\R}$.
\item[\rm{(2)}] If $\sum_{m = 1}^{\infty} \rho_m < \infty$, then $R_* \in (0, 1)$, $\F_{\R} = \F_{\eta}^*$ and
\[
 \E_{\R}(u, v) = \frac 1{R_*}\E^*(u, v) + \frac 1{1 - R_*}\D^I_{\eta}(u, v)
\]
for any $u, v \in \F_{\R}$.
\end{itemize}
\endthm

The idea to prove this theorem will be to show that the restriction of $(\E_{\R}, \F_{\R})$ to $\F_{\eta}$ in the case $\sum_{m=1}^\infty\rho_m = \infty$, respectively $\F_{\eta}^*$ in the case $\sum_{m=1}^\infty\rho_m < \infty$, is again completely symmetric and derived from the same matching pair as $(\E_{\R}, \F_{\R})$. Developing this strategy requires some effort and consists in several steps shown in the subsequent lemmas.\par
We start with two remarks.
\remark
\begin{itemize}[leftmargin=.25in]
\item[\rm{(i)}] By Lemma~\ref{SRF.lemma40}-(1), 
\begin{equation}\label{DRF.eq20}
\sum_{m = 1}^{\infty} \Big(\frac 53\Big)^{m - 1}\c_m = 1 - \prod_{m = 1}^{\infty} (1 - \rho_m)
\end{equation}
and therefore
\[
\sum_{m = 1}^{\infty} \Big(\frac 53\Big)^{m - 1}\eta_m = 1.
\]
\item[\rm{(ii)}] By Proposition~\ref{SRF.prop10},
\[
1 - R_* = r_0(\E_{\R}^I, \F_{\R}^I) = \rho_0.
\]
\end{itemize}
\endremark

\definition\label{DRF.def10}
For each $\R = \{(r_m, \rho_m)\}_{m \ge 1} \in \MP^{\BbN}$ and each $n\geq 0$, define $\R^{(n)} = \{(r_{m + n}, \rho_{m + n})\}_{m \ge 1}$, $R_*^{(n)} = \prod_{m=1}^\infty (1 - \rho_{m + n})$,
\[
\eta^{(n)}_m = \frac{r_{n + 1}{\cdots}r_{n + m - 1}\rho_{n + m}}{1 - R_*^{(n)}}
\]
for $m \ge 1$, and $\eta^{(n)} = \{\eta_m^{(n)}\}_{m \ge 1}$. Moreover, for each $n\geq 0$, define
\[
\F^{(n)} = \begin{cases}
\F_{\eta^{(n)}}&\text{ if }\sum_{m = 1}^{\infty} \rho_m = \infty,\\
\F_{\eta^{(n)}}^*&\text{if }\sum_{m = 1}^{\infty} \rho_m < \infty,
\end{cases}
\]
with $\F_{\eta^{(n)}}$ and $\F^*_{\eta^{(n)}}$ as in Definition~\ref{RFF.def15}.
\enddefinition
\lemma\label{DRF.lemma10}
Let $\R = \{(r_m, \rho_m)\}_{m \ge 1} \in \MP^{\BbN}$ and define for each $n\geq 0$ $\E^{(n)}_{p_1}(u, v) = \E_{\R^{(n)}}(u, v) + u(p_1)v(p_1)$ for any $u, v \in \F_{\R^{(n)}}$. Further, recall the domains $\tF_{\infty}$ and $\F_{\infty}^*$ introduced in Definition~\ref{RFF.def16}.

\begin{itemize}[leftmargin=.25in]
\item[\rm{(1)}] If $\sum_{m = 1}^{\infty}  \rho_m = \infty$, $\F^{(n)}$ is the closure of $\tF_{\infty}$ with respect to the inner product $\E^{(n)}_{p_1}$. 
\item[\rm{(2)}] If $\sum_{m = 1}^{\infty}  \rho_m < \infty$, $\F^{(n)}$ is the closure of $\F_{\infty}^*$ with respect to the inner product $\E^{(n)}_{p_1}$.
\end{itemize}
In either case, $\F^{(n)} \subseteq \F_{\R^{(n)}}$.

\endlemma

\demo
(1)~It suffices to show the case $n = 0$. Let us assume first that $\sum_{m = 1}^{\infty} \rho_m = \infty$. Then, $R_* = 0$ and $\E_{\R} = \E_{\R}^I = \D_{\eta}^I$. Consider now $u \in \F_\eta$, i.e. $u \in \tF$, $\D_{\eta}^I(u, u) <~\infty$ and there exists $\{u_n\}_{n \ge 1} \subseteq \tF_{\infty}$ such that $\lim_{n \to \infty}\D_{\eta}^I(u - u_n, u - u_n) = 0$ and $\lim_{n \to \infty}u_n(x) = u(x)$ for any $x \in K$. Then, $\{u_n\}_{n \ge 1}$ is a Cauchy sequence in $(\F_{\R}, \E_{p_1})$. Since $(\D_{\eta}^I, \F_{\R})$ is a resistance form, there exists $\tilde{u} \in \F_\R$ such that $\E_{p_1}(\tilde{u} - u_n, \tilde{u} - u_n) \to 0$ and $u_n(x) \to \tilde{u}(x)$ as $n \to \infty$ for any $x \in K$. Therefore, $u = \tilde{u} \in \F_{\R}$ and hence it belongs to the closure of $\tF_{\infty}$ with respect to the inner product $\E_{p_1}$. Conversely, it is easy to see that the closure of $\tF_{\infty}$ with respect to $\E_{p_1}$ is a subset of $\F_{\eta}$. Thus, $\F_{\eta}$ is the closure of $\tF_{\infty}$ with respect to the inner product $\E_{p_1}$ and in particular $\F_{\eta} \subseteq \F_{\R}$. If $\sum_{m = 1}^{\infty} \rho_m < \infty$, it follows from Theorem~\ref{CRF.thm10} that
\begin{equation}\label{DRF.eq70}
\E_{\R}(u, v) = \frac{1}{R_*}\E^*(u, v) + \frac 1{1 - R_*}\D_{\eta}^I(u, v)
\end{equation}
for any $u, v \in \F_{\R}$. Consider now $u\in \F^*_\eta$, i.e. $u \in \tF \cap \F^{\SS}$, $\D_{\eta}^I(u, u) < \infty$ and there exists $\{u_n\}_{n \ge 1} \subseteq \F_{\infty}^*$ such that $\lim_{n \to \infty}\E^*(u - u_n, u - u_n) = \lim_{m \to \infty} \D_{\eta}^I(u - u_n, u - u_n) = 0$ and $\lim_{n \to \infty}u_n(x) = u(x)$ for any $x \in K$. Similar arguments as the previous case imply that $u$ belongs to $\F_{\R}$ and $\E_{p_1}(u - u_n, u - u_n) \to 0$ as $n \to \infty$, hence $\F^{(0)}$ is a subset of the closure of $\F_{\infty}^*$ with respect to $\E_{p_1}$. The converse inclusion is straightforward and the desired statement follows.
\enddemo

\definition\label{DRF.def20}
For each $\R = \{(r_m, \rho_m)\}_{m \ge 1} \in \MP^{\BbN}$ and any $n \ge 1$, define $\E^{(n)} = \E_{\R^{(n)}}|_{\F^{(n)} \times \F^{(n)}}$.
\enddefinition

\lemma\label{DRF.lemma20}
Let $\R = \{(r_m, \rho_m)\}_{m \ge 1} \in \MP^{\BbN}$. For any $n \ge 0$ and any $i \in S$,
\[
\{u{\circ}G_i~|~u \in \F^{(n)}\} = \F^{(n + 1)}.
\]
\begin{multline*}
\F^{(n)} = \{u~|~u \in C(K),~u{\circ}G_i \in \F^{(n + 1)}\text{ for any }i \in S,\\
~ u|_{e_{ij}} \in H^1(e_{ij})\text{ for any }(i, j) \in B\},
\end{multline*}
and
\[
\E^{(n)}(u, v) = \sum_{i \in S} \frac 1{r_{n + 1}}\E^{(n + 1)}(u{\circ}G_i, v{\circ}G_i) + \frac 1{\rho_{n + 1}}\D^I_1(u, v)
 \]
for any $u, v \in \F^{(n)}$.
\endlemma
\begin{proof}
From Theorem~\ref{thm20} we know that $\{((\delta_n)^{-1}\E_{\R^{(n)}},\F_{\R^{(n)}},\gamma_n)\}_{n\geq 0}$, where $\delta_n=r_1\cdots r_n$ and $\gamma_n=\delta_{n-1}\rho_n$, is the resolution of $(\E_\R,\F_\R)$. This directly implies the last equality of the lemma because $\F^{(n)}\subseteq\F_{\R^{(n)}}$ by Lemma~\ref{DRF.lemma10}. In view of that equality, if $u\in\F^{(n)}$, then $u|_{e_{ij}} \in H^1(e_{ij})$ for any $(i, j) \in B$. In addition, $u\in\F^{(n)}$ implies the existence of a sequence $\{u_k\}_{k\geq 1}$ that approximates $u$, see Definition~\ref{RFF.def15}, so that $\{u_k{\circ} G_i\}_{k\geq 1}$ approximates $u{\circ}G_i$ in the corresponding way and hence $u{\circ }G_i\in\F^{(n+1)}$. 
On the other hand, consider $u\in C(K)$ such that $u{\circ}G_i\in\F^{(n+1)}$ for any $i\in S$ and $u|_{e_{ij}}\in H^1(e_{ij})$ for all $(i,j)\in B$. Our aim is to prove that $u\in\F^{(n)}$. Since $u{\circ}G_i\in\F^{(n+1)}$, $\E^{(n+1)}(u{\circ}G_i,u{\circ}G_i)<\infty$ and by Lemma~\ref{DRF.lemma10} there exists $\{u_{k,i}\}_{k\geq 1}\subseteq\tF_\infty$ (resp. $\F_\infty^*$)  such that 
\begin{equation*}
\lim_{k \to \infty}\E^{(n+1)}(u{\circ}G_i-u_{k,i},u{\circ}G_i-u_{k,i}) = 0\quad\text{and}\quad \lim_{k \to \infty} u_{k,i}(x) = u{\circ}G_i(x).
\end{equation*}
for any $x \in G_i(K)$. For each $k\geq 1$ define $v_k : K\to\BbR$ by
\begin{equation*}
v_k(x):=\left\{\begin{array}{ll}
u_{k,i}{\circ}G_i^{-1}(x)&\text{if }x\in G_i(K),\\
u(x) + \vp^{ij}_k(x) &\text{if }x\in e_{ij}, (i,j)\in B,
\end{array}\right.
\end{equation*}
where $\vp^{ij}_k$ is an affine function on $e_{ij}$ chosen so that $v_k \in C(K)$. Since $\lim_{k \to \infty} \vp_k^{ij}(p_{ij}) = \lim_{k \to \infty} \vp_k^{ij}(p_{ji}) = 0$,  we have $\D_{e_{ij}}(\vp_k^{ij}, \vp_k^{ij}) \to 0$. By construction, $v_k\in C(K)$ and $v_k\in\tF_\infty$ (resp. $\F_\infty^*$) for any $k\geq 1$. Furthermore, $\D_1^I(u-v_n,u-v_n)\to 0$ as $n\to\infty$ and hence $\E^{(n)}(u-v_k,u-v_k)\to 0$ as $k\to\infty$. Moreover, $\lim_{k\to\infty}v_k(p_1)=\lim_{k\to\infty}v_k(G_1(p_1))=\lim_{k\to\infty}u_{k,1}(p_1)=u(p_1)$ and therefore $u\in\F^{(n)}$.\par
It remains to prove that  $\{u{\circ}G_i~|~u \in \F^{(n)}\} = \F^{(n + 1)}$. On the one hand, it follows from the previous discussion that if $u\in\F^{(n)}$, then $u{\circ}G_i\in\F^{(n+1)}$. On the other hand, consider $u\in\F^{(n+1)}$. By Lemma~\ref{DRF.lemma10}, $\F^{(n+1)}\subseteq\F_{\R^{(n+1)}}=\{v{\circ} G_i~|~v\in\F_{\R^{(n)}}\}$ and we can pick $v\in\F_{\R^{(n)}}$ such that $v{\circ}G_i=u$ for any $i\in S$. In particular, $v\in C(K)$, $v{\circ}G_i\in\F^{(n+1)}$ and $v_{|e_{ij}}\in H^1(e_{ij})$ for all $(i,j)\in B$, so that $v\in\F^{(n)}$.
\end{proof}

\lemma\label{DRF.lemma30}
Let $\R = \{(r_m, \rho_m)\}_{m \ge 1} \in \MP^{\BbN}$. Then, $(\E^{(0)}, \F^{(0)})\in\RF_S$ and its resolution is $\{((\delta_m)^{-1}\E^{(m)}, \F^{(m)}, \c_m)\}_{m \ge 0}$, where $\delta_m = r_1{\cdots}r_m$ and $\c_m = \delta_{m - 1}\rho_m$.
\endlemma

\demo
We start by showing that $(\E^{(0)}, \F^{(0)})$ is a resistance form. Condition (RF1) is obvious. Condition (RF2) follows immediately from Lemma~\ref{DRF.lemma10}. Moreover, since $\tF_{\infty}$ already has the property (RF3) and $\tF_{\infty} \subseteq \F^{(0)}$, (RF3) is also fulfilled. Condition (RF4) holds because $\F_{\eta} \subseteq \F_{\R}$, and 
\begin{equation}\label{DRF.eq80}
\sup\bigg\{\frac{|u(x) - u(y)|^2}{\E(u, u)}~\big|~u \in \F_{\eta},~\E_{\R}(u, u) \neq 0\bigg\} \le R(x, y)
\end{equation}
for any $x,y\in K$. \par
It remains to prove (RF5). Suppose first that $\sum_{m= 1}^\infty \rho_m = \infty$. Obviously, $\tF_{\infty}$ has the Markov property. Now, let $\mu$ be a Borel regular probability measure on $K$ that satisfies $\mu(O) > 0$ for any non-empty open set $O$ and $\mu(A) = 0$ for any finite set $A$. Define
\[
\E_{\mu}(u, v) = \E_{\R}(u, v) + \int_{K} |u(x)|^2\mu(dx)
\]
for any $u, v \in \F_{\R}$. Due to the fact that
\[
|u(x) - u(p_1)|^2 \le \E_{\R}(u, u)R(x, p_1) \le C\E_{\R}(u, u),
\]
where $C = \sup_{x \in K} R(x, p_1)$, we can find $C' > 0$ such that
\[
\frac 1{C'}\E_{p_1}(u, u) \le \E_{\mu}(u, u) \le C'\E_{p_1}(u, u)
\]
for any $u \in \F_{\R}$. Therefore, by Lemma~\ref{DRF.lemma10}, $\F^{(0)}$ is the closure of $\tF_{\infty}$ with respect to $\E_{\mu}$ and~\cite[Theorem~3.1.1]{FOT94} implies that $(\E^{(0)}, \F^{(0)})$ is a Dirichlet form on $L^2(K, \mu)$. In particular, $\F^{(0)}$ has the Markov property and hence (RF5) holds in this case. Suppose now that $\sum_{m = 1}^{\infty} \rho_m < \infty$. Replacing $\tF_{\infty}$ by $\F_{\infty}^*$, the previous arguments show that (RF5) holds again. Thus, $(\E^{(0)}, \F^{(0)})$ is a resistance form. \par
Let $R^{(0)}$ be the resistance metric on $K$ associated with $(\E^{(0)}, \F^{(0)})$ that equals the left-hand side of~\eqref{DRF.eq80}. In view of~\eqref{DRF.eq80}, the identity map from $(K, R)$ to $(K, R^{(0)})$ is continuous and since $(K, R)$ is homeomorphic to $(K, d_E)$, it is compact. Therefore, the identity map from $(K, R)$ to $(K, R^{(0)})$ is a homeomorphism. The rest of the statement follows immediately from Lemma~\ref{DRF.lemma20}.
\enddemo

We finally show Theorem~\ref{DRF.thm10} by making use of these preliminary lemmas to prove that any completely symmetric resistance form $(\E_\R,\F_\R)$ actually coincides with the resistance form $(\E^{(0)},\F^{(0)})$ introduced in Definition~\ref{DRF.def10}. The representation of $\E^{(0)}$ as linear combination of $\E^*$ and $\D_\eta^I$ appears in the proof of Lemma~\ref{DRF.lemma10}, while the domain $\F^{(0)}$ is explicitly given in Definition~\ref{RFF.def15}.

\demo[Proof of Theorem~\ref{DRF.thm10}]
Set $\xi_m = r_0(\E^{(m)}, \F^{(m)})$. Then $r_0((\delta_m)^{-1}\E^{(m)}, \F^{(m)}) = \delta_m\xi_m$. By Lemma~\ref{DRF.lemma30}, the resolution of $(\E^{(0)}, \F^{(0)})$ is $\{((\delta_m)^{-1}\E^{(m)}, \F^{(m)}, \c_m)\}_{m \ge 0}$ and the results in Section~\ref{IRF}, in particular Definition~\ref{IRF.def20} and Theorem~\ref{IRF.thm20}, yield
\begin{equation}\label{DRF.eq90}
\R_{(\E^{(0)}, \F^{(0)})} = \Big\{\Big(\frac{\delta_m\xi_m}{\delta_{m - 1}\xi_{m - 1}}, \frac{\c_m}{\delta_{m - 1}\xi_{m - 1}}\Big)\Big\}_{m \ge 1} = \Big\{\Big(r_m\frac{\xi_m}{\xi_{m - 1}}, \frac{\rho_m}{\xi_{m - 1}}\Big)\Big\}_{m \ge 1}.
\end{equation}
Thus, for any $m \ge 1$,
\begin{equation}\label{DRF.eq100}
\frac 53r_m\frac{\xi_m}{\xi_{m - 1}} + \frac{\rho_m}{\xi_{m - 1}} = 1.
\end{equation}
Since $r_m = \frac 35(1 - \rho_m)$, \eqref{DRF.eq100} yields
\begin{equation}\label{DRF.eq110}
(1 - \xi_m)(1 - \rho_m) = 1 - \xi_{m - 1}
\end{equation}
for any $m \ge 1$, and therefore
\begin{equation}\label{DRF.eq120}
\xi_m = \frac{\xi_0 - 1}{(1 - \rho_1)\cdots(1 - \rho_m)} + 1
\end{equation}
for any $m \ge 1$. Now, it suffices to show that  $\xi_m=1$ for any $m\geq0$.\\
Case 1: Assume that $\sum_{m = 1}^{\infty} \rho_m = \infty$. Since $\xi_m > 0$, we have
\begin{equation}\label{DRF.eq130}
1 - \xi_0 < (1 - \rho_1)\cdots(1 - \rho_m)
\end{equation}
for any $m \ge 1$. The limit of the right-hand side of \eqref{DRF.eq130} as $m \to \infty$ is $0$, hence $\xi_0 \ge 1$. On the other hand, it follows from~\eqref{DRF.eq80} that $\xi_0 = r_0(\E^{(0)}, \F^{(0)}) \le r_0(\E_{\R}, \F_{\R}) = 1$. Therefore, $\xi_0 = 1$ and \eqref{DRF.eq120} implies that $\xi_m = 1$ for any $m \ge 1$. Thus, $\R_{(\E^{(0)}, \F^{(0)})} = \R$. Moreover, $r_0(\E^{(0)}, \F^{(0)}) = \xi_0 = 1 = r_0(\E_{\R}, \F_{\R})$ and Corollary~\ref{IRF.cor10} yields $(\E_{\R}, \F_{\R}) = (\E^{(0)}, \F^{(0)})$.\\
Case 2: Assume that $\sum_{m = 1}^\infty \rho_m < \infty$. By \eqref{DRF.eq70} we have that
\[
\E_{\R_{(\E^{(0)}, \F^{(0)})}}(u, v) = \xi_0\E^{(0)}(u, v) = \frac{\xi_0}{R_*}\E^*(u, v) + \frac{\xi_0}{1 - R_*}\D_{\eta}^I(u, v)
\]
for any $u, v \in \F^{(0)}$. Since $\F^* \subseteq \F^{(0)}$, Theorem~\ref{CRF.thm10} and \eqref{DRF.eq90} yield
\[
\sum_{m = 1}^{\infty} \frac{\rho_m}{\xi_{m - 1}} < \infty
\]
as well as
\begin{equation}\label{DRF.eq140}
\frac{R_*}{\xi_0} = \prod_{m = 1}^{\infty} \Big(1 - \frac{\rho_m}{\xi_{m - 1}}\Big).
\end{equation}
On the one hand, in view of~\eqref{DRF.eq120}, we have that $\{\xi_m\}_{m \ge 1}$ converges as $m \to \infty$. Set $\xi = \lim_{m \to \infty} \xi_m$. Now,~\eqref{DRF.eq110} leads to
\[
1 - \frac{\rho_m}{\xi_{m - 1}} = \frac{\xi_m}{\xi_{m - 1}}(1 - \rho_m),
\]
hence by~\eqref{DRF.eq140}, 
\[
\frac{R_*}{\xi_0} = \frac{R_*}{\xi_0}\xi
\]
and therefore $\xi = 1$. On the other hand, it follows from~\eqref{DRF.eq120} that
\[
\xi = \frac{\xi_0 - 1}{R_*} + 1.
\]
This implies $\xi_0 = 1$ and thus $\R_{(\E^{(0)}, \F^{(0)})} = \R$, which shows $(\E^{(0)}, \F^{(0)}) = (\E_{\R}, \F_{\R})$.
\enddemo

The final step to prove Theorem~\ref{RFF.thm00} consists in showing that any (positive) linear combination of $\E^*$ and $\D_{\eta}^I$ can be realized as a completely symmetric resistance form on SSG.

\demo[Proof of Theorem~\ref{RFF.thm00}]
(1)~If $(\E, \F) \in \RF_S$, then Theorem~\ref{IRF.thm20} implies $(\E, \F) = (c\E_{\R}, \F_{\R})$ for some $c > 0$ and $\R \in \MP^{\BbN}$. By Theorem~\ref{DRF.thm10}, there exist $a \ge 0, b > 0$ and a sequence $\eta = \{\eta_m\}_{m \ge 1} \subseteq (0, \infty)$ such that $\eta$ satisfies~\eqref{RFF.eq30} and  $\E(u,v) = a\E^*(u,v) + b\D_{\eta}^I(u,v)$ for any $u,v\in\F$ with $\F$ as in~\eqref{RFF.eq31}.\par
Conversely, let $\eta = \{\eta_m\}_{m \ge 1} \subseteq (0, 1)$ satisfy \eqref{RFF.eq30}. Inductively we may construct a sequence $\{\s_m\}_{m \ge 1}$ such that
\begin{multline*}
\eta_m = \Big(\frac 35\Big)^{m - 1}(1 - \s_1){\cdots}(1 - \s_{m - 1})\s_m\\
= \Big(\frac 35\Big)^{m - 1}\big((1 - \s_1){\cdots}(1 - \s_{m - 1}) - (1 - \s_1){\cdots}(1 - \s_m)\big)
\end{multline*}
for any $m \ge 1$. In view of~\eqref{RFF.eq30}, it follows that $\prod_{m = 1}^{\infty} (1 - \s_m) = 0$, hence $\sum_{m = 1}^{\infty} \s_m = \infty$. Defining $s_m = \frac 35(1 - \s_m)$ for any $m \ge 1$, $\R_* = \{(s_m, \s_m)\}_{m \ge 1} \in \MP^{\BbN}$ and Theorem~\ref{DRF.thm10} yields
\[
\E_{\R_*}(u, v) = \D_{\eta}^I(u, v)
\]
for any  $u, v \in \F_{\R_*} = \F_{\eta}$. Thus, for any $b > 0$, $(b\D_\eta^I, \F_\eta) = (b\E_{\R_*}, \F_{\R_*}) \in \RF_S$ and the case $a = 0$ of Theorem~\ref{RFF.thm00}-(1) is proven.\par
In order to prove the case $a> 0$, choose $\rho_0 \in (0, 1)$ arbitrarily and define $\rho_m$ for $m\ge 1$ by~\eqref{SRF.eq135}. Then, $\rho_m\in (0, 1)$ for any $m\ge 1$. Taking $r_m = \frac 35(1 - \rho_m)$, we have that $\R = \{(r_m, \rho_m)\}_{m \ge 1}\in \MP^{\BbN}$. Now, set $A_m = \rho_0 \prod_{i = 1}^m (1 - \s_i) + (1 - \rho_0)$ and notice that~\eqref{SRF.eq135} leads to
\[
1 - \rho_m = \frac{A_m}{A_{m - 1}},
\]
hence
\[
R_* = \prod_{m = 1}^{\infty} (1 - \rho_m) = \lim_{m \to \infty} A_m = 1 - \rho_0 > 0.
\]
Moreover, by~\eqref{SRF.eq110},
\[
r_1{\cdots}r_{m - 1}\rho_m = \Big(\frac 35\Big)^{m - 1}(1 - \rho_1){\cdots}(1 - \rho_{m - 1})\rho_m = \Big(\frac 35\Big)^{m - 1}\rho_0\s_m\prod_{i = 1}^{m - 1}(1 - \s_i) = \eta_m\rho_0
\]
and Theorem~\ref{DRF.thm10} yields
\[
\E_{\R}(u, v) = \frac 1{1 - \rho_0}\E^*(u, v) + \frac 1{\rho_0}\D_{\eta}^I(u, v)
\]
for any $u, v \in \F_{\R} = \F_{\eta}^*$. Since $\rho_0 \in (0, 1)$ is arbitrary, for every pair $(a, b) \in (0, \infty) \times (0, \infty)$ in the statement of Theorem~\ref{RFF.thm00}-(1), we find $(a\E^*+b\D^I_\eta, \F_\eta^*)=(c\E_\R, \F_\R) \in \RF_S$ by setting $\rho_0 = a/(a + b)$ and $c = ab/(a + b)$.\\
(2)\,\,Let $\eta = \{\eta_m\}_{m \ge 1}$ satisfy~\eqref{RFF.eq30}. Choose any $\rho_0 \in (0, 1)$ and construct $\R_*$ and $\R$ as in (1). Then, it follows that $\L(\R) = \R_*$. Note that $\F_{\eta} = \F_{\R_*}$ and $\F_{\eta}^* = \F_{\R}$, hence Theorem~\ref{SRF.thm20} yields $\F_{\eta} = \F_{\R_*} = \F_{\R}^I = \{u~|~u \in \F_{\eta}^*, \E_{\R}^{\SS}(u, u) = 0\}$. Since $\E_{\R}^{\SS} = \frac 1{R_*}\E^*$, we finally obtain (2).
\enddemo
\bibliographystyle{amsplain}
\bibliography{RefsDeformedSG}

\end{document}